
\font\Bbb=msbm10 
\def\matC{\hbox{\Bbb C}}
\def\matR{\hbox{\Bbb R}}

\def\matN{\hbox{\Bbb N}}

\font\Gros=cmbx10 scaled\magstep1

\magnification \magstep1 

\null  
\bigskip
\bigskip

\centerline{\Gros Points fixes d'applications holomorphes } 
\smallskip
\centerline{\Gros d'un domaine born\'e de ${\matC}^{n}$ dans lui-m\^eme}  
\medskip
\centerline{\Gros Jean-Pierre Vigu\'e} 
\bigskip
{\bf1. Introduction} 
\medskip

Le but de cet article est de regrouper un certain nombre de r\'esultats
 sur les points fixes d'une application holomorphe $ f:D\longrightarrow
 D $ d'un domaine born\'e $ D $ de $ {\matC}^{n} $ (ou plus g\'en\'eralement
 d'une vari\'et\'e hyperbolique) dans lui-m\^eme. Sur cette question,
 nous allons montrer les deux r\'esultats suivants. 
\medskip

{\bf Th\'eor\`eme 1.1.} - {\it Soit \/}$ D $ {\it un domaine born\'e
 convexe de \/}${\matC}^{n}${\it. Soit \/}$ f:D\longrightarrow
 D $ {\it une application holomorphe. Alors l'ensemble \/}  ${\rm Fix} f = \{z\in
 D|f(z)=z\} $ {\it est une sous-vari\'et\'e analytique complexe de \/}
$D${\it, et si \/}  ${\rm Fix} f ${\it est non vide, il existe une r\'etraction
 holomorphe de \/}$ D $ {\it sur \/} ${\rm Fix} f ${\it. En particulier, ceci
 entra\^{\i}ne que \/} ${\rm Fix} f $ {\it est connexe.\/} 
\medskip

La d\'emonstration de ce r\'esultat utilise un r\'esultat de C. Earle
 et R. Hamilton [8] : soit $ D $ est un domaine born\'e de $ {\matC}^{n}
 $ et soit $ f:D\longrightarrow D $ est une application holomorphe
 telle que $ f(D) $ soit relativement compact dans $ D$. Alors il existe
 une constante $ k<1 $ telle que $ f $ soit $ k$-contractante pour la distance
 int\'egr\'ee de Carath\'eodory $ c^{i}_{D}$. Par suite, elle admet
 un point fixe unique $ a\in D$, et $ a={\rm lim}_{n\rightarrow\infty}f^{n}(z_{0})$,
 pour un point $ z_{o} $ quelconque dans $ D$. 
\medskip

Pour montrer ce r\'esultat, nous serons amen\'es \`a d\'efinir la
 m\'etrique infinit\'esimale et la distance de Carath\'eodory et \`a
 montrer un certain nombre de leurs propri\'et\'es. 
\medskip

Dans le cas d'un domaine born\'e $ D $ de $ {\matC}^{n} $ quelconque,
 les choses sont un peu diff\'erentes et nous montrerons le th\'eor\`eme
 suivant. 
\medskip

{\bf Th\'eor\`eme 1.2. }- {\it Soit \/}$ D $ {\it un domaine born\'e
 de \/}${\matC}^{n} $ {\it(ou une vari\'et\'e hyperbolique). 
Soit \/}$ f:D\longrightarrow D $ {\it une application holomorphe. Alors
 l'ensemble \/} ${\rm Fix} f ${\it des points fixes de \/}$ f $ {\it est une
 sous-vari\'et\'e analytique complexe de \/}$ D${\it.\/} 
\medskip

Pour d\'emontrer ce r\'esultat, nous montrerons d'abord l'existence
 d'une r\'etraction holomorphe $ \rho $ sur une sous-vari\'et\'e analytique
 complexe $ V $ de $ D $ telle que $ \rho(D) $ contienne l'ensemble  ${\rm Fix} f $
 des points fixes de $ f $ et que $ f|_{\rho(D)} $ soit un automorphisme
 analytique de $ \rho(D)$. Ensuite, dans le cas d'une vari\'et\'e hyperbolique
 $ X$, nous \'etudierons le groupe des automorphismes analytiques de
 $ X$. Nous montrerons un th\'eor\`eme de lin\'earisation locale du
 groupe d'isotropie d'un point. Ce th\'eor\`eme, qui a son int\'er\^et
 propre, permet d'achever la d\'emonstration du th\'eor\`eme 1.2. Nous
 remarquerons aussi que la sous-vari\'et\'e  ${\rm Fix} f $ n'est pas connexe
 en g\'en\'eral et que ses composantes connexes n'ont pas toutes la
 m\^eme dimension. 
\medskip

Pour conclure cet article, nous montrerons un r\'esultat plus pr\'ecis
 dans le cas o\`u $ D $ est la boule-unit\'e ouverte de $ {\matC}^{n}
 $ pour une norme et est strictement convexe. 
\medskip

Pour les d\'efinitions
 et les propri\'et\'es des fonctions holomorphes sur un ouvert $ U $ de
 $ {\matC}^{n}$, nous renvoyons le lecteur aux livres de L. et B. Kaup [15]
 et de R. Gunning et H. Rossi [10]. Nous allons commencer par montrer le lemme
 de Schwarz et le lemme de Schwarz-Pick. 
\medskip

{\bf2. Le lemme de Schwarz} 
\medskip

Soit $ \Delta=\{z\in{\matC}| |z|<1\} $ le disque-unit\'e ouvert dans
 $ {\matC}$. 
\medskip

{\bf Th\'eor\`eme 2.1. }(Lemme de Schwarz) - {\it Soit \/}
$ f:\Delta\longrightarrow\Delta $ {\it une application holomorphe telle
 que \/}$ f(0)=0${\it. Alors,\/} 

(i) $ \forall z\in\Delta, |f(z)|\leq|z| $ {\it et \/}$|f' (0)|\leq1${\it.\/}

(ii) {\it De plus, s'il existe \/}$ z_{0}\neq0 $ {\it tel que \/}$|f(z_{0})|=|z_{0}|
 $ {\it ou, si \/}$|f'(0)|=1${\it, alors, il existe un nombre complexe
 \/}$\lambda\in{\matC}, |\lambda|=1 $ {\it tel que \/}$ f(z)=\lambda
 z${\it.\/} 
\medskip

{\it Id\'ee de la d\'emonstration. - \/}Comme $ f(0)=0$, on peut consid\'erer
 l'application holomorphe $ g $ d\'efinie par $ g(z)=f(z)/z$. Si on suppose
 que $ |z|=r$, on trouve que $ |g(z)|\leq1/r$, et, d'apr\`es le principe
 du maximum, c'est aussi vrai pour tout $ |z|\leq r$. En faisant tendre
 $ r $ vers $ 1$, on en d\'eduit que, pour tout $ z\in\Delta, |g(z)|\leq1$.
 Ceci montre le (i). 
\medskip

Pour le (ii), s'il existe $ z_{0}\neq0 $ tel que $ |f(z_{o})|=|z_{0}|
 $ (resp. si $ |f'(0)|=1$), on en d\'eduit que $ |g(z_{0})|=1 $ (resp. $|g(0)|=1$).
 On applique encore une fois le pincipe du maximum pour montrer que
 $ g $ est constante \'egale \`a $ \lambda, |\lambda|=1$. Le (ii) s'en
 d\'eduit. 
\medskip

Le lemme de Schwarz permet de trouver quels sont les automorphismes
 analytiques de $ \Delta $ laissant l'origine fixe. 
\medskip

{\bf Corollaire 2.2. }- {\it Soit \/}$ f:\Delta\longrightarrow\Delta
 $ {\it un automorphisme analytique de \/}$\Delta $ {\it tel que \/}$ f(0)=0.
 $ {\it Alors il existe un nombre complexe \/}$\lambda\in{\matC}, |\lambda|=1
 $ {\it tel que \/}$ f(z)=\lambda z${\it.\/} 
\medskip

{\it D\'emonstration. - \/}Remarquons d'abord que $ z\mapsto\lambda
 z $ avec $ |\lambda|=1 $ est bien un automorphisme analytique de $ \Delta$.
 R\'eciproquement, soit $ f $ un automorphisme analytique de $ \Delta
 $ tel que $ f(0)=0$. Soit $ g=f^{-1}$. D'apr\`es le lemme de Schwarz,
 on a, pour tout $ z\in\Delta$, 
$$|g(f(z))|\leq|f(z)|\leq|z|.$$
 Comme $ g(f(z))=z$, on en d\'eduit que $ |f(z)|=|z|$, et, d'apr\`es
 le lemme de Schwarz, $ f(z)=\lambda z$, avec $ |\lambda|=1$. 
\medskip

Ce r\'esultat permet de montrer quels sont les automorphismes analytiques
 du disque-unit\'e ouvert $ \Delta$. Soit $ \Phi_{a} $ l'application holomorphe
 d\'efinie sur $ {\matC}\backslash\{1/\overline {a}\} $ par $ \Phi_{a}(z)=(z-a)/(1-\overline{a}z),
 a\in\Delta$. On dit que $ \Phi_{a} $ est une transformation de M\"obius.
 Supposons que $ a\in\Delta$. Le premier r\'esultat est que la restriction
 de $ \Phi_{a} $ au disque-unit\'e $ \Delta $ est un automorphisme analytique
 de $ \Delta$. En effet, si on suppose que $ |z|=|e^{i\theta}|=1$, on
 a : 
$$|\Phi_{a}(e^{i\theta})|=|e^{i\theta}-a|/|1-\overline{a}e^{i\theta}|=|e^{i\theta}||1-ae^{-i\theta}|/|1-\overline{a}e^{i\theta}|=1.$$
 Comme $ \Phi_{a} $ est visiblement non constante, on en d\'eduit du
 principe du maximum que, pour tout $ z\in\Delta, |\Phi_{a}(z)|<1$.
 Ainsi donc, $ \Phi_{a} $ envoie $ \Delta $ dans $ \Delta$. Pour montrer
 que $ \Phi_{a} $ est un automorphisme analytique de $ \Delta$, il suffit
 de v\'erifier que $ \Phi^{-1}_{a}=\Phi_{-a}$, ce qui est tout \`a fait
 \'el\'ementaire. 
\medskip

On d\'eduit de ces consid\'erations la forme des automorphismes analytiques
 du disque-unit\'e $ \Delta$. 
\medskip

{\bf Th\'eor\`eme 2.3. }- {\it Les automorphismes analytiques du disque-unit\'e
 \/}$\Delta $ {\it sont les applications de la forme :\/} 
$$z\mapsto\lambda\Phi_{a}(z)=\lambda(z-a)/(1-\overline{a}z),
 |\lambda|=1, a\in\Delta\hbox{\rm.}$$

\medskip

 {\it D\'emonstration. - \/}Il est clair que toutes les applications
 d\'efinies ci-dessus sont des automorphismes analytiques du disque-unit\'e
 $ \Delta$. Soit maintenant $ g $ un automorphisme analytique de $ \Delta
 $ et soit $ a=-g^{-1}(0)$. On a : 
$$g(\Phi_{a}(0))=g(g^{-1}(0))=0.$$
 Par suite, $ g{\circ}\Phi_{a}(z)=\lambda z$, avec $ |\lambda|=1$. On
 en d\'eduit que 
$$g(z)=\lambda\Phi_{a}^{-1}(z)=\lambda\Phi_{-a}(z),$$
 et le th\'eor\`eme est d\'emontr\'e. 
\medskip

{\bf3. Le lemme de Schwarz-Pick} 
\medskip

En utilisant la forme explicite des automorphismes analytiques de
 $ \Delta$, on peut g\'en\'e\-ra\-li\-ser le lemme de Schwarz \`a des applications
 holomorphes de $ \Delta $ dans $ \Delta$, sans supposer que l'origine
 est fixe. Plus pr\'ecis\'ement, nous avons le lemme suivant (lemme
 de Schwarz-Pick, voir par exemple S. Dineen[7]). 
\medskip

{\bf Proposition 3.1 }{\it- Soit f une application holomorphe de \/}$\Delta
 $ {\it dans \/}$\Delta${\it. Alors,\/} 

(i) $ \forall z\in\Delta, \forall w\in\Delta,$ 
$$\Big| \,{f(z)-f(w)   \over 1-\overline{f(w)}f(z) }\,\Big| \leq \Big|\,{z-w
      \over 1-\overline{w}z     }\,\Big|\hbox{\rm;}$$

 (ii) $ \forall z\in\Delta,$ 
$$\,{|f'(z)|     \over 1-|f(z)|^{2}  }\, \leq\,{1       \over 1-|z|^{2}
     }\,.$$
 {\it Si \/}$ f\in{\rm Aut }(\Delta)${\it, \/}(i) {\it et \/}(ii) {\it
 sont des \'egalit\'es. De plus, si on a \'egalit\'e dans \/}(i) {\it
 pour un couple de points distincts \/}$ (z,w) $ {\it ou dans \/}(ii)
 {\it pour un \/}$ z\in\Delta${\it, alors \/}$ f $ {\it est un automorphisme
 analytique de \/}$\Delta${\it.\/} 
\medskip

{\it D\'emonstration. - \/}Soit $ w $ fix\'e $ \in\Delta$. On consid\`ere
 $ g=\Phi_{f(w)}{\circ}f{\circ}\Phi_{-w}$. Alors $ g$ est une 
application holomorphe de $\Delta $ dans $\Delta$ 
et $ g(0)=0$. D'apr\`es le lemme de Schwarz, pour tout $ \zeta\in\Delta$,
 $ |g(\zeta)|\leq|\zeta|$. Par suite, pour tout $ z\in\Delta$, 
$$|g(\Phi_{w}(z))|\leq|\Phi_{w}(z)|,$$
 $$|\Phi_{f(w)}{\circ}f{\circ}\Phi_{-w}{\circ}\Phi_{w}(z)|\leq|\Phi_{w}(z)|,$$
 $$|\Phi_{f(w)}(f(z))|\leq|\Phi_{w}(z)|,$$
 ce qui montre le r\'esultat. Si on a \'egalit\'e dans (i) pour deux
 points distincts, on a alors $ |g(\zeta)|=|\zeta|$, pour un certain
 $ \zeta $ non nul, ce qui, d'apr\`es le lemme de Schwarz montre que
 $ g$, et par suite $ f $ sont des automorphismes analytiques de $ \Delta$.
 \medskip

La d\'emonstration du (ii) est tout \`a fait semblable. 
\medskip

{\bf4. La m\'etrique de Poincar\'e} 
\medskip

Les r\'esultats du paragraphe pr\'ec\'edent s'expriment de mani\`ere
 plus agr\'eable en uti\-li\-sant la m\'etrique infinit\'esimale et la
 distance de Poincar\'e. 
\medskip

{\bf D\'efinition 4.1. }- Soient $ z\in\Delta $ et $ v\in{\matC}$. On
 d\'efinit la m\'etrique infinit\'esimale de Poincar\'e 
$$\alpha(z,v)=\,{|v|     \over 1-|z|^{2}   }\,.$$
 Il est tout \`a fait clair que $ \alpha $ est une m\'etrique hermitienne.
 La m\'etrique $ \alpha $ est invariante dans le sens 
suivant : si $ f:\Delta\longrightarrow\Delta $ est une application holomorphe, alors,
 pour tout $ z\in\Delta $ et tout $ v\in{\matC}$, on a : 
$$\alpha(f(z),f'(z).v)\leq\alpha(z,v),$$
 avec \'egalit\'e si $ f $ est un automorphisme analytique de $ \Delta$.
 Ce r\'esultat est une simple traduction du lemme de Schwarz-Pick.
 \medskip

A partir de cette m\'etrique, on peut calculer la longueur d'une courbe
 $ \gamma:[a,b]\longrightarrow\Delta $ de classe $ {\cal
 C}^{1} $ par morceaux 
$$L_{\alpha}(\gamma)= \int^{b}_{a}\alpha(\gamma(t),\gamma'(t)) dt.$$
 On d\'efinit alors une distance int\'egr\'ee $ \omega $ appel\'ee distance
 de Poincar\'e 
$$\omega(z,w)={\rm inf} L_{\alpha}(\gamma)$$
 pour l'ensemble des courbes $ \gamma $ d'origine $ z $ et d'extr\'emit\'e
 $ w$. 
\medskip

Il est facile de v\'erifier que $ \omega $ est une distance sur $ \Delta$.
 C'est une distance invariante dans le sens suivant : 
si $ f:\Delta\longrightarrow\Delta $ est une application holomorphe, alors,
 pour tout $ z $ et pour tout $ w $ appartenant \`a $ \Delta$, on a : 
$$\omega(f(z),f(w))\leq\omega(z,w),$$
 avec \'egalit\'e si $ f $ est un automorphisme analytique de $ \Delta$.
 \medskip

On peut, bien s\^ur, calculer la valeur de $ \omega $ et on trouve.
\medskip

{\bf Proposition 4.2. }{\it- Pour tous \/}$ z $ {\it et \/}$ w $ {\it
 appartenant \`a \/}$\Delta${\it, on a :\/} 
$$\omega(z,w)= \hbox{\rm\ Arg th }\Big|\,{z-w      \over 1-\overline{w}z
     }\,\Big|\hbox{\rm.}$$
\medskip
 
{\it Id\'ee de la d\'emonstration. - \/}Dans la mesure o\`u les automorphismes
 analytiques de $ \Delta $ sont des isom\'etries, il suffit en fait de
 d\'emontrer ce r\'esultat quand $ w=0 $ et quand $ z $ est un nombre r\'eel
 $ 0<z<1$. Soit $ \gamma(t)=\gamma_{1}(t)+i\gamma_{2}(t) $ un chemin de
 classe {\it C\/}$^{1} $ par morceaux d'origine $ 0 $ et d'extr\'emit\'e
 $ z$. On a : 
$$L_{\alpha}(\gamma)=\int_{a}^{b}\,{|\gamma'_{1} (t)+i\gamma'_{2}
 (t)|        \over 1-|\gamma_{1}(t)+i\gamma_{2}(t)|^{2}   }\,dt$$
 $$\hbox{\rm et, comme }\,{|\gamma'_{1} (t)+i\gamma'_{2} (t)|    
    \over 1-|\gamma_{1}(t)+i\gamma_{2}(t)|^{2}   }\,\geq\,{|\gamma'_{1}
 (t)|               \over 1-|\gamma_{1}(t)|^{2}          }\,\hbox{\rm,
 on trouve que}$$
 $$L_{\alpha}(\gamma)\geq\int_{a}^{b}\,{|\gamma'_{1} (t)|      \over
 1-|\gamma_{1}(t)|^{2}  }\, dt$$
 $$\geq\int_{a}^{b}\,{\gamma'_{1} (t)       \over 1- \gamma_{1}(t)^{2}
   }\, dt=\hbox{\rm Arg th }z.$$
 D'autre part, le chemin $ \varepsilon:[0,z]\longrightarrow\Delta
 $ d\'efinie par $ \varepsilon(t)=t $ a exactement pour longueur 
$\hbox{\rm Arg th }z$. Le r\'esultat est d\'emontr\'e. 
\medskip

{\bf5. D\'efinition des pseudom\'etriques et des pseudodistances invariantes}
\medskip

Nous consid\'ererons un domaine $ D $ de $ {\matC}^{n}$. [Moyennant
 quelques changements \'evidents, la m\^eme construction peut \^etre
 faite pour une vari\'et\'e complexe quelconque (voir S. Ko\-baya\-shi
 [16])]. Il y a deux m\'ethodes pour d\'efinir les pseudodistances invariantes
 sur un domaine $ D $ : une en utilisant les applications holomorphes
 de $ D $ dans le disque-unit\'e $ \Delta$, l'autre en utilisant les applications
 holomorphes de $ \Delta $ dans $ D$ (voir [9], [13] et [14]).
 Commen\c cons par les applications
 holomorphes de $ D $ dans $ \Delta$. 
\medskip

On d\'efinit la pseudodistance de Carath\'eodory $ c_{D} $ de la fa\c
 con suivante : pour tous $ z $ et $ w $ dans $ D$, 
$$c_{D}(z,w)=\hbox{\rm sup}_{\varphi\in H(D,\Delta)} \omega(\varphi(z),\varphi(w)),$$
 o\`u $ H(D,\Delta) $ d\'esigne l'ensemble des applications de $ D $ dans
 le disque-unit\'e $ \Delta$. 
\medskip

On dit qu'une application $ d:D\times D\longrightarrow{\matR}^{+}$ 
est une pseudodistance si $ d $ v\'erifie les trois propri\'et\'es
 suivantes : 

(i) $ x=y\Rightarrow d(x,y)=0 ;$ 

(ii) $ \forall x\in D,\forall y\in D, d(x,y)=d(y,x) ;$ 

(iii) $ \forall x\in D,\forall y\in D,\forall z\in D, d(x,z)\leq d(x,y)+d(y,z).$

En fait, $ d $ v\'erifie toutes les propri\'et\'es d'une distance sauf
 peut-\^etre $ d(x,y)=0\Rightarrow x=y.$ 
\medskip

On a le th\'eor\`eme suivant. 
\medskip

{\bf Th\'eor\`eme 5.1. }{\it- \/}$ c_{D} $ {\it est une pseudodistance
 sur \/}$ D${\it. Elle est invariante dans le sens suivant : pour toute
 application holomorphe \/}$ f:D\longrightarrow D' $ {\it,
 pour tout \/}$ x\in D${\it, pour tout \/}$ y\in D${\it, \/}$ c_{D'}(f(x),f(y))\leq
 c_{D}(x,y). $ {\it Ceci entra\^{\i}ne en particulier que les isomorphismes
 analytiques sont des isom\'etries.\/} 
\medskip

La d\'emonstration se fait en remarquant que l'application 
$ \varphi\mapsto \varphi{\circ}f $ envoie $ H(D' ,\Delta) $ dans
 $ H(D,\Delta)$. 
\medskip

On peut d\'efinir maintenant la fonction de Kobayashi $ \delta_{D}$.
 Pour tous $ z $ et $ w $ appartenant \`a $ D$, on d\'efinit 
$$\delta_{D}(z,w) = \hbox{\rm\ inf }\{\omega(\zeta,\eta) |\; \exists\varphi:\Delta
\longrightarrow D \hbox{\rm\ holomorphe telle que }\varphi(\zeta)=z,\varphi(\eta)=w\}.$$
 Pour les m\^emes raisons, $ \delta_{D} $ est une fonction invariante
 au sens pr\'ec\'edent, mais $ \delta_{D} $ n'est pas une pseudodistance,
 car $ \delta_{D} $ ne v\'erifie pas l'in\'egalit\'e triangulaire en
 g\'en\'eral. 
\medskip

On d\'efinit alors la pseudodistance de Kobayashi $ k_{D} $ comme la
 plus grande pseudodistance inf\'erieure ou \'egale \`a $ \delta_{D}$,
 et on montre que 
$$k_{D}(z,w)= \hbox{\rm\ inf }[\delta_{D}(z_{1},z_{2})+ ... +\delta_{D}(z_{n-1},z_{n})],$$
 pour toutes les suites de points $ (z_{1}, ... ,z_{n}) $ telles que
 $ z_{1}=z, z_{n}=w$. 
\medskip

On montre alors que la pseudodistance de Kobayashi $ k_{D} $ est une
 pseudodistance invariante. 
\medskip

Maintenant, nous pouvons faire la m\^eme construction pour les m\'etriques
 infini\-t\'e\-si\-ma\-les. On d\'efinit alors la pseudom\'etrique infinit\'esimale
 de Carath\'eodory (ou de Reiffen-Cara\-th\'eo\-dory) de la fa\c con suivante
 : $ \forall z\in D, \forall v\in{\matC}^{n}, $  
$$E_{D}(z,v)=\hbox{\rm sup }_{\varphi\in H(\Delta,D)} \alpha(\varphi(z),\varphi'
 (z).v)=\hbox{\rm sup }|\varphi' (z).v|.$$
 (Dans le cas d'une vari\'et\'e $ V$, il faut supposer que $ v $ appartient
 \`a l'espace tangent $ T_{z}(V)$). On v\'erifie facilement que $ E_{D}
 $ est une pseudom\'etrique invariante dans le sens suivant : si $ f:D\longrightarrow D' $ 
est une application holomorphe, on a,
 $ \forall z\in D, \forall v\in{\matC}^{n},$ 
$$E_{D'}(f(z),f'(z).v)\leq E_{D}(z,v).$$
 En particulier, les isomorphismes analytiques sont des isom\'etries.
 
\medskip

De mani\`ere duale, on d\'efinit la pseudom\'etrique infinit\'esimale
 de Kobayashi (ou de Royden-Kobayashi): $ \forall z\in D, \forall v\in{\matC}^{n},$ 
$$F_{D}(z,v)=\hbox{\rm inf}_{\varphi\in H(\Delta,D),\varphi(0)=z}\{|\lambda|
 \hbox{\rm\ tel que }\lambda\varphi' (0)=v\}.$$
 On v\'erifie que $ F_{D} $ est elle aussi une pseudom\'etrique invariante.
 \medskip

Pour la pseudom\'etrique infinit\'esimale de Carath\'eodory, on montre
 la proposition suivante. 
\medskip

{\bf Proposition 5.2. }{\it- L'application \/}
$ E_{D}:D\times{\matC}^{n}\longrightarrow{\matR}^{+} $ {\it est continue.\/}
 \medskip

A partir de la pseudom\'etrique infinit\'esimale de Carath\'eodory,
 on d\'efinit la longueur d'un chemin $ \gamma $ de classe {\it C\/}$^{\hbox{\it1\/}}
 $ par morceaux par la formule 
$$L_{E}(\gamma)= \int^{b}_{a}E_{D}(\gamma(t),\gamma'(t)) dt,$$
 et on d\'efinit la pseudodistance int\'egr\'ee de Carath\'eodory
 de la fa\c con suivante: \'etant donn\'es $ z $ et $ w $ appartenant \`a
 $ D$, on d\'efinit $ c^{i}_{D}(z,w) $ comme la borne inf\'erieure des
 longueurs $ L_{E}(\gamma) $ des chemins $ \gamma $ de classe {\it C\/}$^{\hbox{\it1\/}}
 $ par morceaux d'origine $ z $ et d'extr\'emit\'e $ w$. On v\'erifie facilement
 que $ c^{i}_{D} $ est une pseudodistance invariante. 
\medskip

La pseudom\'etrique infinit\'esimale de Kobayashi nous sera moins
 utile dans ces notes. Pour l'\'etudier, il faut d'abord montrer la
 proposition suivante. 
\medskip

{\bf Proposition 5.3. }{\it- L'application \/}$ F_{D}:D\times{\matC}^{n}\longrightarrow{\matR}^{+} $ {\it est semicontinue
 sup\'e\-rieu\-re\-ment.\/} 
\medskip

On peut alors d\'efinir la longueur $ L_{F}(\gamma) $ d'un chemin de
 classe {\it C\/}$^{\hbox{\it1\/}} $ par morceaux puis d\'efinir une
 pseudodistance int\'egr\'ee $ k^{i}_{D}$. En fait, on montre que $ k^{i}_{D}=k_{D}$.
 \medskip

En utilisant le lemme de Schwarz, on montre facilement le th\'eor\`eme
 suivant. 
\medskip

{\bf Th\'eor\`eme 5.4. }{\it- Pour tout \/}$ z\in D${\it, pour tout
 \/}$ w\in D${\it, on a \/} 
$$c_{D}(z,w)\leq c^{i}_{D}(z,w)\leq k_{D}(z,w)\leq\delta_{D}(z,w).$$
 {\it Ce sont toutes des applications invariantes.\/} 
\medskip

Nous allons donner maintenant quelques exemples de calcul de pseudodistances
 invariantes. 

{\bf Exemple 5.5. }- Pour le disque-unit\'e $ \Delta$, on trouve 
$$c_{\Delta}=k_{\Delta}=\omega, E_{\Delta}=F_{\Delta}=\alpha.$$
 Ce r\'esultat, tout \`a fait \'el\'ementaire est une cons\'equence
 du lemme de Schwarz. 
\medskip

{\bf Exemple 5.6. }- Soit $ D={\matC} $ ou $ {\matC}^{n}$. Dans ce
 cas, $ c_{D}, k_{D} $ et $ \delta_{D} $ sont identiquement nulles. De
 m\^eme, $ E_{D} $ et $ F_{D} $ sont identiquement nulles. 
\medskip

Pour montrer que $ \delta_{D} $ (et par suite, $ c_{D} $ et $ k_{D}$) sont
 identiquement nulles, il suffit de remarquer que, $ z $ et $ w $ \'etant
 donn\'es, l'application $ \varphi:\Delta\longrightarrow
 D $ d\'efinie par $ \varphi(\zeta)=z+n\zeta(w-z) $ est telle que $ \varphi(0)=z,
 \varphi(1/n)=w. $ Le r\'esultat s'en d\'eduit. 
\medskip

{\bf Exemple 5.7. }- Soit $ B $ la boule-unit\'e ouverte de $ {\matC}^{n}$ pour
 une norme $ \|$.$\|$. Alors, 

(i) $ c_{B}(0,z)=k_{B}(0,z)=\delta_{B}(0,z)=\omega(0,\|z\|);$ 

(ii) $ E_{B}(0,v)=F_{B}(0,v)=\|v\|.$ 
\medskip

Montrons par exemple (i). Le r\'esultat est trivial pour $ z=0. $ Supposons
 donc $ z\neq0$. On d\'efinit une application holomorphe $ \varphi:\Delta\longrightarrow B $
 par la formule $ \varphi(\zeta)=(\zeta/\|z\|)z.
 $ D'autre part, d'apr\`es le th\'eor\`eme de Hahn-Banach, il existe
 une application $ {\matC}$-lin\'eaire $ \psi $ de norme 1 telle que
 $ \psi(z)=\|z\|$. La restriction de $ \psi $ \`a $ B $ est donc une application
 holomorphe de $ B $ dans $ \Delta$. En utilisant le fait que $ c_{D} $ est
 une pseudodistance invariante, on trouve 
$$c_{\Delta}(\psi(\varphi(0)),\psi(\varphi(z)))\leq c_{B}(\varphi(0),\varphi(\|z\|))\leq
 c_{\Delta}(0,\|z\|).$$
 Comme $ \psi{\circ}\varphi={\rm id}$, on a : 
$$c_{\Delta}(0,\|z\|)\leq c_{B}(0,z)\leq c_{\Delta}(0,\|z\|),$$
 ce qui montre le r\'esultat pour $ c_{D}$. Le reste des r\'esultats
 se montre de mani\`ere semblable et est laiss\'e en exercice. 
\medskip

Ces r\'esultats permettent de calculer les pseudom\'etrques et les
 pseudodistances sur $ B $ lorsque $ B $ est homog\`ene. Ce r\'esultat
 permet aussi pour un domaine born\'e $ D $ de $ {\matC}^{n} $ et, en
 utilisant des boules bien choisies, de montrer le r\'esultat suivant.
 \medskip

{\bf Th\'eor\`eme 5.8. }{\it- Soit \/}$ D $ {\it un domaine born\'e
 de \/}${\matC}^{n}${\it. Alors \/}$ c_{D} $ {\it et \/}$ k_{D} $ {\it
 sont des distances invariantes qui d\'efinissent la topologie de \/}$ D${\it.
 De m\^eme, \/}$ E_{D} $ {\it et \/}$ F_{D} $ {\it sont des m\'etriques
 invariantes.\/} 
\medskip

Rappelons la d\'efinition suivante : on dit que $ G $ est une m\'etrique
 si, pour tout $ x\in D$, il existe un voisinage $ U $ de $ x $ et une constante
 $ k>0 $ telle que, $ \forall y\in U, \forall v\in{\matC}^{n}, G(y,v)\geq
 k\|v\|$. 
\medskip

On dit qu'un domaine (ou plus g\'en\'eralement une vari\'et\'e) $ D
 $ est hyperbolique si $ k_{D} $ est une distance sur $ D$. Alors, d'apr\`es
 un th\'eor\`eme de T. Barth [2], $ k_{D}$, qui est une distance int\'egr\'ee,
 d\'efinit la topologie de $ D$. 
\medskip

Remarquons que les r\'esultats pr\'ec\'edents montrent que tout domaine
 born\'e $ D $ de $ {\matC}^{n} $ est hyperbolique. 
\medskip

Enfin, il nous faut signaler le profond th\'eor\`eme de L. Lempert
 [17 et 18] (voir aussi M. Jarnicki et P. Pflug [14]). 
\medskip

{\bf Th\'eor\`eme 5.9. }{\it- Soit \/}$ D $ {\it un domaine born\'e
 convexe de \/}${\matC}^{n}${\it. Alors, \/}$ c_{D}=k_{D}=\delta_{D}.
 $ {\it De m\^eme, \/}$ E_{D}=F_{D}.$ 
\medskip

{\bf6. Points fixes d'applications holomorphes} 
\medskip

Comme nous l'avons d\'ej\`a dit, les applications holomorphes sont
 contractantes, au sens large, pour la pseudodistance int\'egr\'ee
 de Carath\'eodory. Nous allons montrer maintenant un r\'esultat de
 C. Earle et R. Hamilton [8] qui montre, dans certains cas, que l'application
 est strictement contractante pour la pseudodistance int\'egr\'ee de
 Carath\'eodory. Plus pr\'eci\-s\'e\-ment, nous avons le th\'eor\`eme suivant.
 \medskip

{\bf Th\'eor\`eme 6.1. }{\it- Soit \/}$ D $ {\it un domaine born\'e
 de \/}${\matC}^{n}${\it, et soit \/}$ f:D\longrightarrow
 D $ {\it une application holomorphe telle que \/}$ f(D) $ {\it soit relativement
 compact dans \/}$ D${\it. Alors il existe une constante \/}$ k<1 $ {\it
 telle que \/}$ f $ {\it soit \/}$ k${\it-contractante pour la distance
 int\'egr\'ee de Carath\'eodory \/}$ c^{i}_{D}${\it. Par suite, \/}$ f
 $ {\it admet un unique point fixe \/}$ a\in D${\it, et \/}
$ a={\rm lim}_{n\rightarrow\infty} f^{n}(z_{0})${\it, o\`u \/}$ z_{0} $ {\it
 est un point quelconque de \/}$ D${\it.\/} 
\medskip

{\it D\'emonstration. \/}- Du fait que $ f(D)\subset\subset D$, on
 d\'eduit l'existence de deux nombres r\'eels $ r $ et $ R $ strictement
 positifs tels que , $ \forall z\in D$, 
$$B(f(z),r)\subset D\subset B(f(z),R).$$
 Soit $ z\in D $ fix\'e et soit $ t=r/R$. L'application holomorphe 
$$\zeta\mapsto g(\zeta)=f(\zeta)+t(f(\zeta)-f(z))$$
 est une application holomorphe de $ D $ dans $ D$. D'apr\`es les propri\'et\'es
 de la m\'etrique de Carath\'eodory, on a donc : 
$$E_{D}(g(z),g'(z).v)\leq E_{D}(z,v),$$
 ce qui donne 
$$(1+t)E_{D}(f(z),f'(z).v)\leq E_{D}(z,v).$$
 Ainsi donc, $ f $ est $ 1/(1+t)$-contractante pour $ E_{D} $ et $ c^{i}_{D}$.
 
\medskip

On choisit alors un point $ z_{0}\in D$. Si on note $ z_{n}=f^{n}(z_{0})$,
 on a : 
$$c^{i}_{D}(z_{n},z_{n+1})\leq1/(1+t)^{n}c^{i}_{D}(z_{0},z_{1}).$$
 De mani\`ere classique, on en d\'eduit que $ (z_{n})_{n\in{\matN}}
 $ est une suite de Cauchy dans $ D $ et que, pour tout $ n\geq1$, $ z_{n}\in
 \overline{f(D)} $ qui est un compact de $ D$. Sa limite est l'unique
 point fixe de $ f $ dans $ D$. 
\medskip

{\bf7. Lin\'earisation locale d'une r\'etraction holomophe} 
\medskip

Nous aurons besoin du th\'eor\`eme suivant (H. Cartan [6]). 
\medskip

{\bf Th\'eor\`eme 7.1. }{\it- Soit \/}$ U $ {\it un ouvert de \/}
${\matC}^{n} $ {\it et soit \/}$\rho:U\longrightarrow U $ {\it
 une r\'etraction holomorphe, c'est-\`a-dire une application holomorphe
 telle que \/}$\rho^{2}=\rho{\circ}\rho=\rho${\it. Soit \/}$ z_0 $
 {\it un point tel que \/}$\rho(z_{0})=z_{0}${\it. Alors il existe
 une carte locale \/}$ u $ {\it d'un voisinage \/}$ U(z_{0}) $ {\it sur
 un voisinage \/}$ V $ {\it de \/}$ 0 $ {\it tel que \/}$ u(z_{0})=0${\it,
 et que, dans la carte \/}$ u${\it, \/}$\rho $ {\it soit lin\'eaire \'egal
 \`a un projecteur, (ce qui signifie que \/}$ u{\circ}\rho{\circ}u^{-1}
 $ {\it est la restriction \`a \/}$ V $ {\it d'un projecteur lin\'eaire).\/}
 
\medskip

{\it D\'emonstration. - \/}Quitte \`a faire un changement de coordonn\'ees,
 on peut supposer que $ \rho(0)=0$. Soit $ \rho=\sum^{\infty}_{n=1} P_{n}
 $ le d\'eveloppement de $ \rho $ en s\'eries de polynomes homog\`enes.
 Le fait que $ \rho^{2}=\rho $ entra\^{\i}ne que $ P_{1}{\circ}P_{1}=P_{1}$.
 Soit 
$$u={\rm id}+(2P_{1}-{\rm id}){\circ}(\rho-P_{1})\hbox{\rm.}$$
 On a: 
$$u'(0)={\rm id}+(2P_{1}-{\rm id}){\circ}(P_{1}-P_{1})={\rm id}.$$
 Le th\'eor\`eme d'inversion locale montre que $ u $ est bien une carte
 locale au voisinage de $ 0$. Etudions maintenant $ u{\circ}\rho$. On
 a : 
$$u{\circ}\rho=\rho+(2P_{1}-{\rm id}){\circ}(\rho-P_{1}){\circ}\rho$$
 $$=\rho+(2P_{1}-{\rm id}){\circ}(\rho^{2}-P_{1}{\circ}\rho)$$
 $$=\rho+2P_{1}{\circ}\rho-2P_{1}{\circ}\rho-\rho+P_{1}{\circ}\rho=P_{1}{\circ}\rho.$$
 D'autre part, 
$$P_{1}{\circ}u=P_{1}+(2P_{1}-P_{1}){\circ}(\rho-P_{1})$$
 $$=P_{1}+P_{1}{\circ}\rho-P_{1}=P_{1}{\circ}\rho.$$
 On a donc 
$$u{\circ}\rho=P_{1}{\circ}u.$$
 On en d\'eduit 
$$u{\circ}\rho{\circ}u^{-1}=P_{1},$$
 et le th\'eor\`eme est d\'emontr\'e. 
\medskip

On d\'eduit de ce r\'esultat le corollaire suivant. 
\medskip

{\bf Corollaire 7.2. }{\it- Soit \/}$ U $ {\it un ouvert de \/}
${\matC}^{n} $ {\it et soit \/}$\rho:U\longrightarrow U $ {\it
 une r\'etraction holomorphe. Alors l'ensemble \/}$\rho(U) $ {\it(qui
 est \'egal \`a l'ensemble des points fixes de \/}$\rho${\it) est une
 sous-vari\'et\'e analytique complexe de \/}$ U${\it.\/} 
\medskip

Le r\'esultat est \'evident parce que, dans la carte $ u$, $ \rho(U)
 $ est \'egal \`a l'intersection de $ V $ et d'un sous-espace vectoriel
 de $ {\matC}^{n}$. 
\medskip

{\bf8. Points fixes d'applications holomorphes dans un domaine born\'e
 convexe de }${\matC}^n$ 
\medskip

On d\'eduit des consid\'erations pr\'ec\'edentes le th\'eor\`eme suivant
 (P. Mazet et J.-P. Vigu\'e [19]). 
\medskip

{\bf Th\'eor\`eme 8.1. }{\it- Soit \/}$ D $ {\it un domaine born\'e
 convexe de \/}${\matC}^{n}${\it. Soit \/}$ f:D\longrightarrow
 D $ {\it une application holomorphe. Alors l'ensemble \/} ${\rm Fix} f = \{z\in
 D|f(z)=z\} $ {\it des points fixes de \/}$ f $ {\it est une sous-vari\'et\'e
 analytique complexe connexe de \/}$ D${\it, et si \/} ${\rm Fix} f $ {\it
 est non vide, il existe une r\'etraction holomorphe \/}$\rho:D\longrightarrow{\rm Fix}  f${\it.\/} 
\medskip

{\it D\'emonstration. \/}- On peut bien s\^ur supposer que  ${\rm Fix} f $
 est non vide. Pour tout nombre r\'eel $ \lambda, 0\leq\lambda<1$,
 et pour tout $ a\in D$, soit 
$$z\mapsto f_{\lambda,a}(z)=a+\lambda(f(z)-a).$$
 Soit $ r>0 $ tel que $ B(a,r) $ soit contenu dans $ D$. En \'ecrivant

$$f_{\lambda,a}(z)=(1-\lambda)a+\lambda f(z),$$
 on montre que $ D $ contient la boule de centre $ f_{\lambda,a}(z) $ et
 de rayon $ (1-\lambda)r$. Ainsi, $ f_{\lambda,a}(D) $ est relativement
 compact dans $ D$. Par suite, $ f_{\lambda,a} $ admet un point fixe unique
 $ \varphi_{\lambda}(a) $ et on a : 
$$\varphi_{\lambda}(a)=\hbox{\rm lim}_{n\rightarrow\infty}f^{n}_{\lambda,a}(b),$$
 o\`u $ b $ est un point quelconque de $ D$. D'autre part,
 $ a\mapsto f^{n}_{\lambda,a}(b) $ est holomorphe, et, d'apr\`es
 le th\'eor\`eme 6.1, on peut v\'erifier que ces applications convergent
 uniform\'ement sur tout compact vers $ \varphi_{\lambda} $ qui est donc
 holomorphe. Remarquons qu'en fait, on peut le v\'erifier directement
 en appliquant le th\'eor\`eme de Montel : on peut trouver $ \psi $ tel
 que $ f^{n_{j}}_{\lambda,a}(b) $ tende vers $ \psi $ uniform\'ement sur
 tout compact. Mais $ \psi=\varphi_{\lambda}$,
 et tout s'en d\'eduit facilement.
 \medskip

On fait alors tendre $ \lambda $ vers $ 1$. D'apr\`es le th\'eor\`eme
 de Montel, il existe une suite $ \lambda_{n}\rightarrow 1
 $ telle que $ \varphi_{\lambda_{n}} $ tende vers $ \varphi$, o\`u $ \varphi
 $ est une application holomorphe de $ D $ dans $ \overline{D}$. Remarquons
 d'abord que, si $ a $ est un point fixe de $ f$, $ f_{\lambda,a}(a)=a$.
 Par suite, $ \varphi_{\lambda,a}(a)=a $ et $ \varphi(a)=a$. Comme $ D
 $ est un domaine convexe born\'e et est, par suite, taut, le fait que
 $ f $ admet au moins un point fixe entra\^{\i}ne que $ \varphi $ est une
 application holomorphe de $ D $ dans $ D$. Pour montrer que $ \varphi $ est
 une r\'etraction holomorphe de $ D $ sur l'ensemble des points fixes
 de $ f$, il suffit de montrer que, pour tout $ b\in D, \varphi(b) $ est
 un point fixe de $ f$, c'est-\`a-dire que $ f(\varphi(b))=\varphi(b)$.
 Or, on sait que 
$$f_{\lambda,b}(\varphi_{\lambda,b}(b))=\varphi_{\lambda,b}(b),$$
 $$(1-\lambda)b+\lambda f(\varphi_{\lambda,b}(b))=\varphi_{\lambda,b}(b).$$
 Faisons tendre $ \lambda_{n} $ vers $ 1$. Alors $ \varphi_{\lambda_{n},b}(b)
 $ tend vers $ \varphi(b)$. Comme $ f $ est continue, on en d\'eduit que
 
$$f(\varphi(b))=\varphi(b).$$
 Ainsi, $ \varphi $ est bien une r\'etraction holomorphe sur l'ensemble
 des points fixes de $ f$, et d'apr\`es le corollaire 7.2,  ${\rm Fix} f $ est
 bien une sous-vari\'et\'e analytique connexe de $ D$. 
\medskip

Maintenant, si on ne suppose pas que $ D $ est un domaine born\'e de
 $ {\matC}^{n}$, alors l'ensemble des points fixes de $ f $ n'est pas,
 en g\'en\'eral, une sous-vari\'et\'e analytique de $ D$. 
\medskip

{\bf Exemple 8.2. }- Consid\'erons le cas de $ {\matC}^{n}$. Soit
 $ f $ l'application holomorphe de $ {\matC}^{n} $ dans $ {\matC}^{n} $ d\'efinie
 par 
$$f(z_{1},...,z_{n})=(z_{1}+f_{1}(z_{1},...,z_{n}),...,z_{n}+f_{n}(z_{1},...,z_{n}))\hbox{\rm,}$$
 o\`u $ f_{1},...,f_{n} $ sont des applications holomorphes
 de $ {\matC}^{n} $ dans $ {\matC}$. Alors l'ensemble des points fixes de $ f $ est
 \'egal \`a l'ensemble des z\'eros communs des $ f_{1},...,f_{n}$. En
 g\'en\'eral, ce n'est pas une sous-vari\'et\'e analytique complexe
 de $ {\matC}^{n}$. 
\medskip

D'autre part, m\^eme dans le cas d'un domaine born\'e $ D $ de $ {\matC}^{n}$, 
l'ensemble  ${\rm Fix} f $ n'est pas, en g\'en\'eral, un ensemble
 connexe, et ses composantes connexes n'ont pas toujours la m\^eme
 dimension. 
\medskip

{\bf Exemple 8.3. }- Soit $ R>1$, et soit $ A=\{z\in{\matC}|1/R<|z|<R\}$.
 Alors l'application $ f:A\longrightarrow A $ d\'efinie
 par $ f(z)=1/z $ admet deux points fixes $ 1 $ et $ -1$. Ce n'est pas un
 ensemble connexe. 
\medskip

{\bf Exemple 8.4. }- Consid\'erons le domaine born\'e suivant 
$$D=\{(x,y)\in{\matC}^{2}|1/2<|x|<2, |xy^{2}|<1\}.$$
 Soit $ f:D\longrightarrow D $ d\'efinie par $ f(x,y)=(1/x,xy)$.
 L'ensemble des points fixes de $ f $ est \'egal \`a $ A_{1}\cup A_{-1}$,
 o\`u 
$$A_{1}=\{(x,y)\in D|x=1\},A_{-1}=\{(x,y)\in D|x=-1,y=0\}.$$
 Il est clair que $ A_{1} $ et $ A_{-1} $ sont de dimension diff\'erentes.
 Cependant, on remarque que, dans ce cas aussi,  ${\rm Fix} f $ est une sous-vari\'et\'e
 analytique complexe de $ D$. Nous allons montrer que, si $ D $ est un
 domaine born\'e de $ {\matC}^{n}$, ou plus g\'en\'eralement une
 vari\'et\'e hyperbolique, et si $ f:D\longrightarrow
 D $ est une application holomorphe, l'ensemble des points fixes de $ f
 $ est une sous-vari\'et\'e analytique complexe de $ D$. 
\medskip

Nous allons montrer ce r\'esultat en utilisant un th\'eor\`eme sur
 la lin\'earisation locale d'un automorphisme analytique de $ D $ au
 voisinage d'un point fixe. Pour commencer, nous aurons besoin de certains
 r\'esultats sur les automorphismes analytiques d'un domaine born\'e.
 
\medskip

{\bf9. Automorphismes analytiques d'un domaine born\'e} 
\medskip

Dans ce paragraphe, nous allons montrer quelques propri\'et\'es du
 groupe des automorphismes analytiques d'un domaine born\'e (voir par 
exemple [20]). 
\medskip

Soit donc $ D $ un domaine born\'e de $ {\matC}^{n}$. On munit l'ensemble
 $ H(D,D) $ des applications holomorphes de $D$ dans $D$ 
de la topologie de la convergence uniforme sur tout compact
 de $ D$. Le groupe ${\rm Aut }(D) $ des automorphismes analytiques de $ D $ qui
 est contenu dans $ H(D,D) $ sera muni de la topologie induite. 
\medskip

Remarquons d'abord que sur $ H(D,D)$, cette topologie co\"{\i}ncide
 avec la topologie de la convergence uniforme sur un compact d'int\'erieur
 non vide. Plus pr\'ecis\'ement, nous avons la proposition suivante.
 
\medskip

{\bf Proposition 9.1. }{\it- Soit \/}$ D $ {\it un domaine born\'e
 de \/}${\matC}^{n}${\it, et soit \/}$ K $ {\it un compact d'int\'erieur
 non vide contenu dans \/}$ D${\it. Alors, sur \/}$ H(D,D)${\it, la topologie
 de la convergence sur tout compact de \/}$ D $ {\it co\"{\i}ncide avec
 la topologie de la convergence uniforme sur \/}$ K${\it.\/} 
\medskip

{\it D\'emonstration. - \/}Pour montrer la proposition 9.1, il suffit
 de d\'emontrer le r\'esultat suivant : soit $ (f_{n})_{n\in{\matN}}
 $ une suite de fonctions holomorphes appartenant \`a $ H(D,D) $ convergeant
 uniform\'ement sur $ K $ vers $ f\in H(D,D)$. Alors $ f_{n} $ converge
 vers $ f $ uniform\'ement sur tout compact de $ D$. D'apr\`es le th\'eor\`eme
 de Montel, $ H(D,\overline{D}) $ est compact. On peut donc extraire
 de la suite $ f_{n} $ une suite $ f_{n_{j}} $ convergeant uniform\'ement
 sur tout compact vers $ g\in H(D,\overline{D})$. Mais $ f $ et $ g $ co\"{\i}ncident
 sur l'int\'erieur de $ K$. D'apr\`es le th\'eor\`eme de prolongement
 analytique, $ f=g$. Ainsi donc, $ f $ est le seul point adh\'erent \`a
 la suite $ f_{n}$. Comme $ H(D,\overline{D}) $ est compact, ceci entra\^{\i}ne
 que $ f_{n} $ converge vers $ f $ uniform\'ement sur tout compact de $ D$.
 \medskip

A l'aide de cette proposition, on montre assez facilement le th\'eor\`eme
 suivant. 
\medskip

{\bf Th\'eor\`eme 9.2. }{\it- \/}(i) {\it L'application \/}$ H(D,\overline{D})\times
 H(D,D)\longrightarrow H(D,\overline{D})${\it, qui, \`a
 \/}$ (f,g) $ {\it associe \/}$ f{\circ}g${\it, est continue ;\/} 

(ii) {\it L'application de \/}${\rm Aut }(D) $ {\it dans lui-m\^eme \/}
$ f\mapsto f^{-1}${\it est continue.\/} 

{\it Par suite, \/}${\rm Aut }(D) $ {\it est un groupe topologique.\/} 
\medskip

{\it D\'emonstration. - \/} (i) Les espaces consid\'er\'es \'etant
 m\'etriques, on peut utiliser des suites. Soient $ f_{n}\rightarrow
 f, g_{n}\rightarrow g. $ Montrons que $ f_{n}{\circ}g_{n}\rightarrow f{\circ}g$. 
\medskip

Choisissons un compact $ K $ d'int\'erieur non vide contenu dans $ D$.
 Comme $ g\in H(D,D)$, $ g(K) $ est compact, et on peut trouver un compact
 $ L $ contenu dans $ D $ tel que, pour un certain $ \delta>0$, $ L $ contienne
 $ g(K)_{\delta}=\{z\in D|d(z,g(K))<\delta\}$. Soit $ z\in K$. On a :
 
$$\|f_{n}(g_{n}(z))-f(g(z))\|\leq\|f_{n}(g_{n}(z))-f(g_{n}(z))\|+\|f(g_{n}(z))-f(g(z))\|.$$
 Pour $ n $ assez grand, $ g_{n}(z)\in L$. On a donc 
$$\|f_{n}(g_{n}(z))-f(g_{n}(z))\|\leq\|f-f_{n}\|_{L}.$$
 D'autre part, $ f $ \'etant continue sur le compact $ g(K) $ est uniform\'ement
 continue sur $ g(K)$. Par suite, pour tout $ \varepsilon>0$, il existe
 $ \eta>0 $ tel que $ \|g_{n}(z)-g(z)\|<\delta $ entra\^{\i}ne $ \|f(g_{n}(z))-f(g(z))\|<\varepsilon.
 $ Ceci d\'emontre que $ f_{n}(g_{n}(z))\rightarrow
 f(g(z)) $ uniform\'ement sur $ K$. D'apr\`es la proposition 9.1, ceci
 entra\^{\i}ne que $ f_{n}{\circ}g_{n} $ converge vers $ f{\circ}g $ uniform\'ement
 sur tout compact de $ D$. 
\medskip

(ii) Soit $ f_{n} $ une suite d'\'el\'ements de ${\rm Aut }(D)$. Soit $ g_{n}=f_{n}^{-1}$,
 et on veut montrer que $ g_{n}$ converge vers $ f^{-1}$. Comme $ g_{n}
 $ est une suite de fonctions holomorphes born\'ees, on sait qu'il existe
 une suite extraite $ g_{n_{j}} $ convergeant vers $ g$. D'apr\`es le
 (i), $ g_{n_{j}}{\circ}f_{n_{j}} $ converge vers $ g{\circ}f$. Comme
 $ g_{n_{j}}$=$f_{n_{j}}^{-1}, g_{n_{j}}{\circ}f_{n_{j}}={\rm id}$. Par suite,
 $ g{\circ}f={\rm id}$. En multipliant \`a droite par $ f^{-1}$, on trouve
 que $ g=f^{-1}$. Ainsi, le seul point adh\'erent \`a la suite $ g_{n}
 $ est $ f^{-1}$. Par suite, $ g_{n} $ converge vers $ f^{-1}$, et le th\'eor\`eme
 est d\'emontr\'e. 
\medskip

Si on consid\`ere une suite $ f_{n} $ d'automorphismes d'un domaine
 born\'e $ D $ convergeant vers une limite $ f$, $ f $ n'est pas forc\'ement
 un automorphisme analytique de $ D$. En effet, dans le cas du disque-unit\'e
 $ \Delta$, il suffit de consid\'erer la suite $ f_{n} $ d'automorphismes
 analytiques d\'efinie par 
$$f_{n}(z)=\,{z+(1-1/n)   \over 1+(1-1/n)z  }\,,$$
 qui converge uniform\'ement sur tout compact vers $ 1$. C'est dans
 un certain sens le seul probl\`eme possible comme le montre le th\'eor\`eme
 suivant. 
\medskip

{\bf Th\'eor\`eme 9.3. }{\it- Soit \/}$ D $ {\it un domaine born\'e
 de \/}${\matC}^{n}${\it, et soit \/}$ f_{n}\in{\rm Aut }(D) $ {\it une suite
 d'automorphismes analytiques de \/}$ D $ {\it convergeant vers \/}$ f\in
 H(D,\overline{D})${\it. S'il existe \/}$ a\in D $ {\it tel que \/}$ f(a)\in
 D${\it, \/}$ f $ {\it est un automorphisme analytique de \/}$ D${\it.\/}
\medskip

Comme pr\'ec\'edemment, on fait la d\'emonstration en consid\'erant
 la suite $ g_{n}=f_{n}^{-1} $ et en montrant que $ g_{n} $ converge vers
 $ f^{-1}$. 
\medskip

{\bf Corollaire 9.4. }{\it- Soit \/}$ D $ {\it un domaine born\'e de
 \/}${\matC}^{n} $ {\it et soit \/}$ a $ {\it un point de \/}$ D${\it.
 Soit \/}$ K $ {\it un compact de \/}$ D${\it. Alors \/}${\rm Aut }_{a,K}(D)=\{f\in{\rm Aut }(D)|f(a)\in
 K\} $ {\it est compact. En particulier, \/} ${\rm Aut }(D) $ {\it est localement
 compact et le groupe \/} ${\rm Aut }_{a}(D) $ {\it d'isotropie du point \/}$ a $ 
 {\it est compact.\/} 
\medskip

{\it D\'emonstration. - \/}Soit $ f_{n} $ une suite d'\'el\'ements de
 ${\rm Aut }_{a,K}(D)$. On peut en extraire une suite $ f_{n_{j}} $ convergeant
 vers $ f\in H(D,\overline{D})$. Comme $ K $ est compact, $ f(a)\in K$.
 D'apr\`es le th\'eor\`eme 9.3, $ f\in{\rm Aut }(D)$. Ainsi, ${\rm Aut }_{a,K}(D)
 $ est compact. 
\medskip

{\bf10. Propri\'et\'es du groupe des automorphismes analytiques d'un
 domaine born\'e} 
\medskip

Nous avons d'abord le th\'eor\`eme d'unicit\'e de H. Cartan [5]. 
\medskip

{\bf Th\'eor\`eme 10.1. }{\it- Soit \/}$ D $ {\it un domaine born\'e
 de \/}${\matC}^{n} $ {\it et soit \/}$ f\in H(D,D) $ {\it tel que \/}$ f(a)=a,
 f' (a)={\rm id}${\it. Alors, \/}$ f={\rm id}$. 
\medskip

{\it D\'emonstration. - \/}Soit $ f(z)=a+\sum^{\infty}_{p=1} P_{p}(z-a)
 $ le d\'eveloppement de $ f $ en s\'eries de polyn\^omes homog\`enes.
 Comme $ f'(a)={\rm id}$, ceci entra\^{\i}ne que $ P_{1}(z-a)=f'(a).(z-a)=(z-a)$.
 Si $ f $ n'est pas \'egal \`a l'identit\'e, on peut \'ecrire 
$$f(z)=a+(z-a)+P_{k}(z-a)+\sum^{\infty}_{p=k+1} P_{p}(z-a),$$
 o\`u $ P_{k} $ est le premier terme non nul du d\'eveloppement d'ordre $ k\geq2$. On v\'erifie
 facilement que, si on it\`ere $ f$, on trouve comme d\'eveloppement
 en s\'eries de polyn\^omes homog\`enes de $ f^n$ 
$$f^{n}(z)=a+(z-a)+nP_{k}(z-a)+...$$
 Les in\'egalit\'es de Cauchy montrent que $ \|nP_{k}\| $ est born\'e
 par une constante $ M$, et est donc nul. Contradiction. Le th\'eor\`eme
 est d\'emontr\'e. 
\medskip

On en d\'eduit le corollaire suivant. 
\medskip

{\bf Corollaire 10.2. }{\it- Soit \/}$ D $ {\it un domaine born\'e de
 \/}${\matC}^{n} $ {\it et soient \/}$ f\in H(D,D) $ {\it et \/}$ g\in{\rm Aut }(D)${\it.
 Si \/}$ f(a)=g(a) $ {\it et \/}$ f' (a)=g' (a) $ {\it alors \/}$ f=g${\it.\/}
 
\medskip

Pour d\'emontrer ce r\'esultat, il suffit d'appliquer le th\'eor\`eme
 10.1 \`a $ g^{-1}{\circ}f$. 
\medskip

On d\'eduit de ces consid\'erations le th\'eor\`eme suivant. 
\medskip

{\bf Th\'eor\`eme 10.3. }{\it- Soit \/}$ D $ {\it un domaine born\'e
 de \/} ${\matC}^{n}${\it. Alors, l'application \/}$\varphi $ {\it
 de \/}${\rm Aut }(D) $ {\it dans \/}$ D\times$GL$({\matC}^{n}) $ {\it d\'efinie
 par \/} 
$$f\mapsto(f(a),f'(a))$$
 {\it est injective et est un hom\'emorphisme de \/}${\rm Aut }(D) $ {\it
 sur son image.\/} 
\medskip

{\it D\'emonstration. - \/}Le fait que $ \varphi $ est injective a \'et\'e
 d\'emontr\'e dans le corollaire 10.2. La continuit\'e de $ \varphi
 $ est imm\'ediate. Montrons maintenant que $ \varphi $ est bicontinue.
 Soit $ f_{n} $ et $ f $ appartenant \`a ${\rm Aut }(D) $ et supposons 
que $ f_{n}(a)\rightarrow f(a), f' _{n}(a)\rightarrow f'(a)$. Montrons que 
$ f_{n}\rightarrow f $ uniform\'ement 
sur tout compact de $ D$.
 On peut extraire de la suite $ f_{n} $ une suite $ f_{n_{j}} $ qui converge
 vers $ g\in H(D,\overline{D})$. Comme $ g(a)\in D$, le th\'eor\`eme
 9.3 montre que $ g\in{\rm Aut }(D)$, et on a : $ f(a)=g(a)$, $ f' (a)=g' (a)$.
 Par suite, $ g=f. $ Ainsi, $ f $ est le seul point adh\'erent \`a la suite
 $ f_{n} $ dans le compact ${\rm Aut }_{a,K}(D)$, o\`u $ K $ est un voisinage compact
 de $f(a)$. Par suite, $ f_{n}\rightarrow f $ et le th\'eor\`eme
 est d\'emontr\'e. 
\medskip

Ainsi donc, si $ D $ est un domaine born\'e de $ {\matC}^{n}$, un \'el\'ement
 $ f $ du groupe ${\rm Aut }_{a}(D) $ d'isotropie du point $ a $ est caract\'eris\'e
 par la valeur de $ f' (a)$. En fait, on a mieux et nous allons montrer
 les deux r\'esultats suivants. 
\medskip

{\bf Th\'eor\`eme 10.4. }{\it- Soit \/}$ X $ {\it une vari\'et\'e analytique
 complexe hyperbolique, et soit \/}$ f:X\longrightarrow
 X $ {\it un automorphisme analytique de \/}$ X${\it, et soit \/}$ a $ {\it
 un point fixe de \/}$ f${\it. Alors il existe une carte locale \/}$\varphi
 $ {\it d\'efinie au voisinage de a, telle que \/}$\varphi(a)=0${\it,
 et que, dans la carte \/}$\varphi${\it, \/}$ f $ {\it soit une application
 lin\'eaire (c'est-\`a-dire que \/}$\varphi{\circ}f{\circ}\varphi^{-1}
 $ {\it soit une application lin\'eaire).\/} 
\medskip

{\it D\'emonstration. - \/}Dans la mesure o\`u la question est locale
 et, quitte \`a remplacer $ X $ par une boule pour la distance de Kobayashi
 de centre $ a $ et de rayon suffisamment petit, on peut supposer $ X
 $ un domaine born\'e de $ {\matC}^{n}$. D\'efinissons 
$$\varphi_{n}(z)=(1/n)\sum_{p=0}^{n-1} f' (a)^{-p}.(f^{p}(z)-a),$$
 $ f^{p} $ d\'esigne la $ p-$i\`eme it\'er\'ee de $ f $ et $ f'(a)^{-p}
 $ la $ p-$i\`eme it\'er\'ee de $ f'(a)^{-1}$. Munissons $ {\matC}^{n}
 $ de la norme $ F_{B}(a,.)$. On sait que $ f' (a) $ est une isom\'etrie
 pour cette norme. D'autre part, la suite $ (f^{p}(z)-a) $ est born\'ee.
 Par suite, $ \varphi_{n} $ est une suite de fonctions holomorphes born\'ees.
 On peut utiliser le th\'eor\`eme de Montel et extraire de la suite
 $ \varphi_{n} $ une suite $ \varphi_{n_{j}} $ convergeant uniform\'ement
 sur tout compact vers une application holomorphe $ \varphi$. On sait
 que 
$$\varphi_{n}'(a)=(1/n)\sum_{p=0}^{n-1} f' (a)^{-p}{\circ}f'(a)^{p}=\hbox{\rm
 id}.$$
 D'apr\`es le th\'eor\`eme de Weierstrass, ceci entra\^{\i}ne que
 $ \varphi'(a)={\rm id}$. Le th\'eor\`eme d'inversion locale montre que $ \varphi
 $ est un isomorphisme analytique d'un voisinage de $ a $ sur un voisinage
 de $ 0$. 
\medskip

Pour d\'emontrer le th\'eor\`eme, il suffit maintenant de montrer
 que $ \varphi{\circ}f{\circ}\varphi^{-1}=f'(a)$, ou, ce qui revient
 au m\^eme, que $ \varphi{\circ}f=f' (a){\circ}\varphi$. 
\medskip

Soit $ z\in D$. Calculons $ \varphi_{n}(f(z))$. 
$$\varphi_{n}(f(z))=(1/n)\sum_{p=0}^{n-1} f' (a)^{-p}.(f^{p+1}(z)-a)$$
 $$=f'(a).(1/n)\sum_{p=0}^{n-1} f' (a)^{-(p+1)}.(f^{p+1}(z)-a)$$
 $$=f'(a).(1/n)\sum_{p=1}^{n} f' (a)^{-p}.(f^{p}(z)-a)$$
 $$=f'(a).(1/n)\sum_{p=0}^{n-1} f' (a)^{-p}.(f^{p}(z)-a)+(1/n)(f'(a)^{-n}.(f^{n}(z)-a)-(z-a))$$
 Quand $ n $ tend vers l'infini, $ (1/n)(f'(a)^{-n}.(f^{n}(z)-a)-(z-a))
 $ tend vers $ 0$. D'autre part, $ f' (a).(1/n_{j})\sum_{p=0}^{n_{j}-1}f'(a)^{-p}.(f^{p}(z)-a)
 $ tend vers $ f'(a).\varphi(z)$. On trouve donc $ f(\varphi(z))=f' (a).\varphi(z)$,
 et le th\'eor\`eme est d\'emontr\'e. 
\medskip

Ce r\'esultat sera suffisant pour d\'emontrer notre th\'eor\`eme sur
 les points fixes. Cependant, nous pouvons montrer un r\'esultat plus
 fort que nous allons \'enoncer maintenant. 
\medskip

{\bf Th\'eor\`eme 10.5. }{\it- Soit \/}$ X $ {\it une vari\'et\'e hyperbolique
 complexe et soit \/}$ a $ {\it un point de \/}$ X${\it. Alors il existe
 une carte locale \/}$\varphi $ {\it d\'efinie au voisinage de a, telle
 que \/}$\varphi(a)=0${\it, et que, dans la carte \/}$\varphi${\it,
 \/}${\rm Aut }_{a}(D) $ {\it soit un sous-groupe du groupe lin\'eaire (c'est-\`a-dire
 que, pour tout \/}$ f\in{\rm Aut }_{a}(D)${\it, \/}$\varphi{\circ}f{\circ}\varphi^{-1}
 $ {\it soit une application lin\'eaire).\/} 
\medskip

{\it D\'emonstration. - \/}Comme pr\'ec\'edemment, on peut supposer
 que $ X $ est un domaine born\'e $ D $ de $ {\matC}^{n}$. Comme ${\rm Aut }_{a}(D)
 $ est compact, le groupe ${\rm Aut }_{a}(D) $ peut \^etre muni d'une mesure
 de Haar $ \mu $ invariante \`a gauche et \`a droite de masse totale
 $ 1 $ (voir N. Bourbaki [4]). On d\'efinit $ \varphi $ par la formule 
$$\varphi(z)=\int_{{\rm Aut}_{a}(D)} g' (a)^{-1}.(g(z)-a)d\mu(z).$$
 Il suffit de d\'eriver sous le signe somme pour v\'erifier que $ \varphi'(a)={\rm id}$.
 Par suite, $ \varphi $ est une carte locale au voisinage de $ a$. 
\medskip

Soit $ f\in{\rm Aut }_{a}(D)$. Calculons $ \varphi{\circ}f $ : 
$$\varphi(f(z))=\int_{{\rm Aut}_{a}(D)} g' (a)^{-1}.(g(f(z))-a)d\mu(g)$$
 $$=f'(a).\int_{{\rm Aut}_{a}(D)} [(g{\circ}f)' (a)]^{-1}.((g{\circ}f)(z)-a)d\mu(g)=f'(a).\varphi(z),$$
 en utilisant l'invariance de $ \mu $ par translation \`a droite. Le
 th\'eor\`eme est d\'emontr\'e. 
\medskip

Dans le cas des domaines cercl\'es born\'es de $ {\matC}^{n}$, nous
 avons un r\'esultat meilleur que nous allons maintenant d\'emontrer.
 Commen\c cons par donner la d\'efiniton suivante. 
\medskip

{\bf D\'efinition 10.6. }- On dit qu'un domaine $ D $ de $ {\matC}^{n}
 $ est cercl\'e si l'origine 0 appartient \`a $ D $ et si, pour tout nombre
 complexe $ \lambda $ de module 1, pour tout $ z\in D, \lambda z\in D$.
 
\medskip

{\bf Th\'eor\`eme 10.7. }{\it- Soient \/}$ D_{1} $ {\it et \/}$ D_{2}
 $ {\it deux domaines cercl\'es de \/}${\matC}^{n} $ {\it et supposons
 que \/}$ D_{1} $ {\it est born\'e. Soit \/}$ f:D_{1}\longrightarrow
 D_{2} $ {\it un isomorphisme analytique tel que \/}$ f(0)=0${\it. Alors,
 \/}$ D_{2} $ {\it est born\'e et \/}$ f $ {\it est la restriction \`a
 \/}$ D_{1} $ {\it d'un automorphisme lin\'eaire de \/}${\matC}^{n}${\it.\/}
 \medskip

{\it D\'emonstration. - \/}Montrons d'abord que, pour tout $ \theta\in{\matR}$, 
pour tout $ z\in D_{1},$ 
$$f(e^{i\theta}z)=e^{i\theta}f(z).$$
 Pour cela, consid\'erons l'application holomorphe
 $ g $ de $ D_{1} $ dans $ D_{1}
 $ d\'efinie par 
$$g(z)=f^{-1}(e^{-i\theta}f(e^{i\theta}z)).$$
 Comme $ D_{1} $ et $ D_{2} $ sont cercl\'es, $ g $ est un automorphisme
 analytique de $ D_{1}$. On a : $ g(0)=0, g' (0)={\rm id}$. Le th\'eor\`eme
 d'unicit\'e de H. Cartan montre que $ g(z)=z$, ce qui entra\^{\i}ne
 que $ f(e^{i\theta}z)=e^{i\theta}f(z).$ 
\medskip

Consid\'erons maintenant le d\'eveloppement de $ f $ en s\'eries de
 polyn\^omes homog\`enes au voisinage de l'origine 
$$f(z)=\sum_{n=1}^{\infty} P_{n}(z),$$
 qui converge uniform\'ement sur tout compact contenu dans $ D_{1}$.
 Soit $ z\in D_{1} $ et consid\'erons le compact r\'eunion des $ e^{i\theta}z$,
 pour $ \theta\in{\matR}$. Il est facile de voir que l'int\'egrale
 
$$(1/2\pi)\int_{0}^{2\pi}f(e^{i\theta}z)e^{-i\theta}d\theta$$
 est \'egale \`a $ P_{1}(z)$. D'un autre c\^ot\'e, en utilisant l'\'egalit\'e
 $ f(e^{i\theta}z)=e^{i\theta}f(z)$, on trouve que cette int\'egrale
 vaut $ f(z)$. Ce r\'esultat entra\^{\i}ne facilement que $ f $ est la
 restriction de l'application lin\'eaire $ P_{1}$. 
\medskip

Ce r\'esultat montre en particulier que l'ensemble des points fixes
 d'un automorphisme analytique d'un domaine cercl\'e born\'e $ D $ de
 $ {\matC}^{n} $ laissant l'origine fixe est l'intersection de $ D $ avec
 un sous-espace vectoriel $ V $ de $ {\matC}^{n}$. Ce sous-espace $ V
 $ est justement le sous-espace propre correspondant \`a la valeur propre
 $ 1 $ de $ f=f' (0)$. 
\medskip

Ce th\'eor\`eme permet aussi de montrer un r\'esultat d\'emontr\'e
 d'abord par H. Poincar\'e: en dimension sup\'erieure ou \'egale \`a
 $ 2$, il existe des domaines born\'es de $ {\matC}^{n} $ simplement
 connexes et non analytiquement isomorphes. 
\medskip

{\bf Th\'eor\`eme 10.8. }{\it- Pour tout \/}$ n\geq2${\it, le polydisque
 \/}$\Delta^{n}=\{(z_{1},...,z_{n})\in{\matC}^{n}||z_{i}|<1,\forall
 i\} $ {\it et la boule hermitienne \/}
$ B_{n}=\{(z_{1},...,z_{n})\in{\matC}^{n}||z_{1}|^{2}+...+|z_{n}|^{2}<1\} $ {\it ne sont pas analytiquement
 isomorphes.\/} 
\medskip

{\it Id\'ee de la d\'emonstration. \/}- En utilisant des transformations
 de M\"obius, il est facile de montrer que le polydisque $ \Delta^{n}
 $ est homog\`ene (c'est aussi vrai pour la boule hermitienne mais moins
 \'evident). S'il existe un isomorphisme $ f $ de $ \Delta^{n} $ sur $ B_{n}$,
 on peut donc supposer que $ f(0)=0$. D'apr\`es le th\'eor\`eme
 10.7, $ f $ serait donc lin\'eaire. Il est facile de voir qu'un tel
 isomorphisme lin\'eaire ne peut pas exister. On peut remarquer par
 exemple que la fronti\`ere $ \partial B_{n} $ de $ B_{n}
 $ est une sous-vari\'et\'e analytique r\'eelle de $ {\matC}^{n} $ alors
 que ce n'est pas le cas pour la fronti\`ere $ \partial\Delta^{n} $ de
 $ \Delta^{n}$. Cependant, $ f $ serait un isomorphime lin\'eaire de $ \partial\Delta^{n}
 $ sur $ \partial B_{n}$, ce qui est impossible. 
\medskip

{\bf11. Suite des it\'er\'ees d'une application holomorphe} 
\medskip

Pour la suite de notre \'etude, nous allons consid\'erer (
comme E. Bedford [3] et M. Abate [1]) la suite
 des it\'er\'ees d'une application holomorphe $ f$. Plus pr\'ecis\'ement,
 on consid\`ere un domaine born\'e $ D $ de $ {\matC}^{n} $ et une application
 holomorphe $ f:D\longrightarrow D$. On consid\`ere alors
 la suite des it\'er\'ees $ f^{n}(=f{\circ}...{\circ}f$ $ n $ fois) de $ f
 $ qui est une suite d'applications holomorphes de $ D $ dans $ D$. Comme
 $ D $ est born\'e, d'apr\`es le th\'eor\`eme de Montel, on peut trouver
 une suite extraite $ (f^{n_{j}}) $ qui converge vers une application
 holomorphe $ g:D\longrightarrow\overline{D}$. 
\medskip

{\bf Th\'eor\`eme 11.1. }{\it- Supposons que la suite \/}$ f^{n} $ {\it
 v\'erifie la propri\'et\'e \/}(H) {\it suivante : pour toute suite
 extraite \/}$ f^{n_{j}} $ {\it convergente vers \/}$ g${\it, \/}$ g $ {\it
 est une application holomorphe de \/}$ D $ {\it dans \/}$ D${\it. Alors
 il existe une r\'etraction holomorphe \/}$\rho:D\longrightarrow
 D $ {\it et une suite extraite \/}$ f^{n_{j}} $ {\it qui converge vers
 \/}$\rho${\it. De plus, \/}$ f|_{\rho(D)} $ {\it est un automorphisme
 analytique de \/}$\rho(D) $ {\it et \/}$\rho(D) $ {\it contient l'ensemble
 des points fixes de \/}$ f${\it.\/} 
\medskip

{\it D\'emonstration. \/}- Commen\c cons par extraire une sous-suite
 $ f^{n_{j}} $ qui converge vers $ F:D\longrightarrow D$.
 On peut supposer, quitte \`a extraire une sous-suite que $ n_{j+1}-2n_{j}
 $ tend vers $ +\infty$. En appliquant plusieurs fois le th\'eor\`eme
 de Montel, on peut trouver une application 
$ \sigma:{\matN}\longrightarrow{\matN} $ strictement croissante
 telle que la suite $ h_{j}=f^{n_{\sigma(j)+1}-n_{\sigma(j)}}
 $ converge vers une application holomorphe $ \rho:D\longrightarrow
 D$. De m\^eme, on suppose que $ f^{n_{\sigma(j)+1}-2n_{\sigma(j)}} $ converge
 vers $ g$. Enfin, supposons que $ f^{n_{\sigma(j)+1}-n_{\sigma(j)}-1}
 $ converge vers $ h$. 
\medskip

En \'ecrivant 
$$f^{n_{\sigma(j)+1}}\hbox{\rm=}f^{n_{\sigma(j)+1}-n_{\sigma(j)}}{\circ}f^{n_{\sigma(j)}}\hbox{\rm=}f^{n_{\sigma(j)}}{\circ}f^{n_{\sigma(j)+1}-n_{\sigma(j)}},$$
 et, en faisant tendre $ j $ vers $ +\infty$, on trouve 
$$F=\rho{\circ}F=F{\circ}\rho.$$
 Soit $ z=F(x) $ un \'el\'ement appartenant \`a l'image de $ F$. On a
 : 
$$\rho(z)=\rho(F(x))=F(x)=z.$$
 Par suite, $ \rho $ est \'egal \`a l'identit\'e sur l'image de $ F$.

Consid\'erons maintenant $ f^{n_{\sigma(j)+1}-n_{\sigma(j)}}$. En \'ecrivant
 
$$f^{n_{\sigma(j)+1}-n_{\sigma(j)}}=f^{n_{\sigma(j)+1}-2n_{\sigma(j)}}{\circ}f^{n_{\sigma(j)}}=f^{n_{\sigma(j)}}{\circ}f^{n_{\sigma(j)+1}-2n_{\sigma(j)}},$$
 et, en passant \`a la limite, on trouve 
$$\rho=g{\circ}F=F{\circ}g.$$
 Ceci prouve que l'image de $ \rho $ est contenue dans l'image de $ F$,
 et on a vu que, sur l'image de $ F$, $ \rho $ est \'egal \`a l'identit\'e.
 Ainsi, $ \rho $ est bien une r\'etraction holomorphe. 
\medskip

Consid\'erons maintenant $ f^{n_{\sigma(j)+1}-n_{\sigma(j)}+1}$. En
 \'ecrivant 
$$f^{n_{\sigma(j)+1}-n_{\sigma(j)}+1}=f^{n_{\sigma(j)+1}-n_{\sigma(j)}}{\circ}f=f{\circ}f^{n_{\sigma(j)+1}-n_{\sigma(j)}},$$
 et, en faisant tendre $ j $ vers $ +\infty$, on trouve 
$$f{\circ}\rho=\rho{\circ}f.$$
 Ceci prouve que l'image de $ \rho $ est stable par $ f$. 
\medskip

Consid\'erons enfin $ f^{n_{\sigma(j)+1}-n_{\sigma(j)}-1} $ qui converge
 vers $ h$. Nous avons 
$$(f^{n_{\sigma(j)+1}-n_{\sigma(j)}})^{2}=f^{n_{\sigma(j)+1}-n_{\sigma(j)}+1}{\circ}f^{n_{\sigma(j)+1}-n_{\sigma(j)}-1}$$
 $$=f^{n_{\sigma(j)+1}-n_{\sigma(j)}-1}{\circ}f^{n_{\sigma(j)+1}-n_{\sigma(j)}+1},$$
 et, en faisant tendre $ j $ vers $ +\infty$, on trouve 
$$\rho^{2}=f{\circ}\rho{\circ}h=h{\circ}\rho{\circ}f=\rho{\circ}h{\circ}f.$$
 Ceci prouve que $ f|_{{\rm Im}\rho} $ est un automorphisme analytique
 de Im $ \rho$. Son inverse est la restriction de $ \rho{\circ}h $ \`a
 l'image de $ \rho$. 
\medskip

Sur la construction de $ f^{n}$, il est clair que, si $ z $ est un point
 fixe de $ f$, alors $ f^{n}(z)=z $ et, par suite, $ \rho(z)=z. $ Ainsi,
 $ z $ appartient \`a l'image de $ \rho$. Le th\'eor\`eme est d\'emontr\'e.
 \medskip

Rappelons qu'une suite $ f_{n} $ d'applications d'une vari\'et\'e
 analytique complexe $ X $ dans une vari\'et\'e complexe $ Y $ est compactement
 divergente si, \'etant donn\'es deux compacts $ K\subset X $ et $ L\subset
 Y,\exists n_{0}\in{\matN} $ tel que, $ \forall n\geq n_{0}, f_{n}(K)\cap
 L=\emptyset$. D'autre part, on dit qu'une vari\'et\'e analytique complexe
 $ X $ est taut (au sens de H. Wu [24]) si, pour toute suite $ f_{n} $ d'applications
 du disque-unit\'e $ \Delta $ dans $ X$, la suite $ f_{n} $ admet une sous-suite
 $ f_{n_{j}} $ convergente dans $ H(\Delta,X) $ ou une sous-suite $ f_{n_{j}}
 $ compactement divergente. 
\medskip

Il est classique qu'une vari\'et\'e analytique complexe hyperbolique
 compl\`ete (c'est-\`a-dire, compl\`ete pour la distance de Kobayashi
 $ k_{X}$) est taut et que toute vari\'et\'e complexe taut est hyperbolique.
 En particulier, il r\'esulte des r\'esultats pr\'ec\'edents la proposition
 suivante. 
\medskip

{\bf Proposition 11.2. }{\it- Soit \/}$ D $ {\it un domaine born\'e
 taut de \/}${\matC}^{n}${\it. Soit \/}$ f $ {\it une application holomorphe
 de \/}$ D $ {\it dans \/}$ D${\it. Supposons qu'il existe \/}$ x\in D
 $ {\it tel que la suite \/}$ f^{n}(x) $ {\it soit relativement compacte
 dans \/}$ D${\it. Alors la suite \/}$ f^{n} $ {\it des it\'er\'ees de
 \/}$ f $ {\it v\'erifie l'hypoth\`ese \/}(H). 
\medskip

En particulier, si on consid\`ere un domaine born\'e $ D $ taut de
 $ {\matC}^{n} $ et une application $ f $ de $ D $ dans $ D $ ayant un point fixe
 dans $ D$, la suite des it\'er\'ees $ f^{n} $ de $ f $ v\'erifie l'hypoth\`ese
 (H). Cependant, ceci n'est pas tout \`a fait suffisant pour l'\'etude
 des points fixes car nous d\'esirons \'etudier ausi le cas o\`u $ D
 $ n'est pas taut. Pour cela, nous avons la proposition suivante. (Nous
 noteros $ B_{k}(a,r) $ la boule de centre $ a $ et de rayon $ r $ pour la
 distance de Kobayashi). 
\medskip

{\bf Proposition 11.3. }{\it- Soit \/}$ D $ {\it un domaine born\'e
 de \/}${\matC}^{n} $ {\it et soit \/}$ f:D\longrightarrow
 D $ {\it une application holomorphe de \/}$ D $ {\it dans \/}$ D $ {\it
 admettant un point fixe \/}$ a${\it. Alors, pour tout \/}$ r>0${\it,
 \/}$ f(B_{k}(a,r))\subset B_{k}(a,r)${\it. De plus, il existe \/}$ r>0
 $ {\it tel que la suite \/}$ (f|_{B_{k}(a,r)})^{n} $ {\it v\'erifie la
 condition \/}(H). 
\medskip

{\it D\'emonstration. \/}- En effet, soit $ r>0 $ tel que la boule $ B_{k}(a,r)
 $ soit isomorphe \`a un domaine born\'e et relativement compacte dans
 $ D$. Soit $ g_{n}=(f|_{B_{k}(a,r)})^{n} $ et soit $ z\in B_{k}(a,r)$.
 Alors $ k_{X}(a,z)=\rho<r. $ Si on consid\`ere une sous-suite $ g_{n_{j}}
 $ convergente, la limite $ g(z) $ appartient \`a l'adh\'erence de $ B_{k}(a,r)
 $ qui est contenue dans X. Par passage \`a la limite, on en d\'eduit
 que $ k_{X}(a,g(z))\leq\rho$, et $ g $ est une application holomorphe
 de $ B_{k}(a,r) $ dans lui-m\^eme. 
\medskip

{\bf12. Points fixes d'une application holomorphe} 
\medskip

Nous pouvons maintenant achever l'\'etude de l'ensemble des points
 d'une application holomorphe $ f $ d'un domaine born\'e (ou plus g\'en\'eralement
 d'une vari\'et\'e hyperbolique) dans lui-m\^eme. Plus pr\'ecisement,
 nous avons le th\'eor\`eme suivant (voir J.-P. Vigu\'e [23]). 
\medskip

{\bf Th\'eor\`eme 12.1. }{\it- Soit \/}$ X $ {\it une vari\'et\'e hyperbolique
 complexe (ou un domaine born\'e de \/}${\matC}^{n}${\it).
 Soit \/}$ f:X\longrightarrow X $ {\it une application holomorphe. Alors,
 l'ensemble \/} ${\rm Fix} f $ {\it des points fixes de \/}$ f $ {\it est une
 sous-vari\'et\'e analytique complexe de \/}$ X${\it. De plus, si \/}$ a
 $ {\it est un point fixe de \/}$ f${\it, la dimension de \/} ${\rm Fix} f $ {\it
 au voisinage de \/}$ a $ {\it est \'egale \`a la dimension du sous-espace
 propre correspondant \`a la valeur propre \/}$ 1 $ {\it de \/}$ f' (a)${\it.\/}
 \medskip

{\it D\'emonstration. - \/}Le r\'esultat est de nature locale. Il
 suffit de d\'emontrer que, \'etant donn\'e $ a\in {\rm Fix} f $, il existe
 un voisinage $ U $ de $ a $ tel que  ${\rm Fix} f \cap U $ soit une sous-vari\'et\'e
 analytique de $ U$. 
\medskip

Soit $ a\in {\rm Fix} f $. Quitte \`a remplacer $ X $ par une boule pour la
 distance de Kobayashi de rayon suffisamment petit, on peut supposer
 que $ X $ est un domaine born\'e  de $ {\matC}^{n}$. On applique alors
 la proposition 11.3, ce qui permet ensuite d'appliquer le th\'eor\`eme
 11.1. Il existe un voisinage born\'e $ V $ de $ a $ dans $ {\matC}^{n}
 $ stable par $ f $ qui v\'erifie les hypoth\`eses du th\'eor\`eme 11.1.
 Par suite, il existe une r\'etraction holomorphe $ \rho $ de $ V $ sur
 une sous-vari\'et\'e $ W $ de $ V$, $ W $ contient les points fixes de
 $ f|_{V} $ et $ f|_{W} $ est un automorphisme analytique de $ W$. Quitte
 \`a diminuer encore une fois la taille de $ W$, on peut supposer que
 $ W $ est isomorphe \`a un domaine born\'e de $ {\matC}^{n} $ et appliquer
 le th\'eor\`eme 10.4. Il existe une carte locale d\'efinie au voisinage
 de $ a $ tel que, dans cette carte, $ f $ soit lin\'eaire \'egal \`a $ f'
 (a). $ Le th\'eor\`eme s'en d\'eduit. 
\medskip

Nous allons nous int\'eresser maintenant aux points fixes d'une application
 holomorphe $ f $ d'une boule-unit\'e de $ {\matC}^{n} $ dans elle-m\^eme
 qui laisse l'origine fixe. Commen\c cons par un principe du maximum
 fort. 
\medskip

{\bf13. Principe du maximum fort} 
\medskip

Si on consid\`ere un domaine born\'e $ U $ de $ {\matC}^{n} $ et une
 fonction holomorphe $ f:U\longrightarrow{\matC}$, on
 a le principe du maximum : si $ |f| $ admet un maximum local en un point
 $ a $ de $ U$, $ f $ est constante. 
\medskip

Maintenant, si on consid\`ere une application holomorphe d'un ouvert
 $ U $ de $ {\matC}^{n} $ dans $ {\matC}^{m}$, muni d'une norme $ \|.\|$,
 le principe du maximum ne s'applique pas en g\'en\'eral, $ \|f\| $ peut
 avoir un maximum local sans que $ f $ soit constante, comme le montre
 l'exemple suivant : soit $ {\matC}^{2} $ muni de la norme 
$$\|(z_{1},z_{2})\|=\hbox{\rm sup(}\|z_{1}\|,\|z_{2}\|).$$
 L'application $ f:\Delta\longrightarrow{\matC}^{2} $ d\'efinie
 par $ f(\zeta)=(1,\zeta) $ v\'erifie $ \|f(z)\|=1, \forall\zeta\in\Delta$,
 et $ f $ n'est pas constante. 
\medskip

Commen\c cons par montrer le lemme suivant (voir par exemple [9] ou [12]). 
\medskip

{\bf Lemme 13.1. }{\it- Soit \/}$f \in H(\Delta,\overline{\Delta})${\it.
 Alors, pour tout \/}$ z\in\Delta${\it,\/} 
$$2|z||f(0)|+(1-|z|)|f(z)-f(0)|\leq2|z|.$$
 {\it D\'emonstration. \/}- S'il existe $ z $ tel que $ |f(z)|=1$, $ f
 $ est constante et la formule est vraie. On peut donc supposer que
 $ f(z)\in\Delta, \forall z\in\Delta. $ On peut alors utiliser le lemme
 de Schwarz-Pick (proposition 3.1). On a, $ \forall z\in\Delta $ : 
$$\Big|\,{f(z)-f(0)   \over 1-\overline{f(0)}f(z) }\,\Big| \leq \Big|z\Big|,$$
 ce qui donne 
$$|f(z)-f(0)|\leq|z||1-\overline{f(0)}f(z)|.$$
 Majorons $ |1-\overline{f(0)}f(z)|$. 
$$|1-\overline{f(0)}f(z)|\leq|1-|f(0)|^{2}|+||f(0)|^{2}-\overline{f(0)}f(z)|$$
 $$\leq|1-|f(0)|^{2}|+|\overline{f(0)}||f(0)-f(z)|$$
 $$\leq|1-|f(0)|^{2}+|f(0)||f(0)-f(z)|.$$
 On en d\'eduit 
$$|f(z)-f(0)|\leq|z|(|1-|f(0)|^{2})+|z||f(0)||f(z)-f(0)|.$$
 Par suite, 
$$(1-|z||f(0)|)|f(z)-f(0)|\leq|z|(1-|f(0)|^{2})$$
 $$\leq|z|(1-|f(0)|)(1+|f(0)|)\leq2|z|(1-|f(0)|).$$
 Ceci donne 
$$2|z||f(0)|+(1-|z||f(0)|)|f(z)-f(0)|\leq2|z|,$$
 et comme $ 1-|z|\leq(1-|z||f(0)|)$, le r\'esultat s'en d\'eduit. 
\medskip

Nous avons maintenant le lemme suivant. 
\medskip

{\bf Lemme 13.2. }{\it- Soit \/}$ B $ {\it la boule-unit\'e ouverte
 de \/}${\matC}^{n} $ {\it pour une norme \/}$\|.\|${\it, et soit \/}$ f
 $ {\it une application holomorphe de \/}$\Delta $ {\it dans \/}${\matC}^{n} $ 
{\it telle que \/}$ f(\Delta) $ {\it soit contenu dans \/}$\overline{B}${\it.
 Alors\/} 
$$\|f(0)+\zeta(f(z)-f(0))\|\leq1$$
 {\it pour tout \/}$ z\in\Delta\backslash\{0\} $ {\it et pour tout \/}$\zeta\in{\matC} $
 {\it tel que \/}$|\zeta|\leq\,{1-|z|    \over 2|z|    }\,.$ 
\medskip

{\it D\'emonstration. - \/}Si $ n=1$, on d\'eduit du lemme 13.1 que
 
$$|f(0)|+\,{1-|z|  \over 2|z|  }\, |f(z)-f(0)|\leq1.$$
 Pour tout $ \zeta\in{\matC} $ {\it tel que \/}$|\zeta|\leq\,{1-|z|
    \over 2|z|    }\,$, on a $ \|f(0)+\zeta(f(z)-f(0)\|\leq1$, ce qui
 est le r\'esultat annonc\'e. 

Dans le cas $ n\geq2$, faisons la d\'emonstration par l'absurde. Supposons
 qu'il existe $z\in\Delta\backslash\{0\} $ 
et $ \zeta\in{\matC}$, $ |\zeta|\leq\,{1-|z|
    \over 2|z|    }\, $ tel que $ |f(0)|+\,{1-|z|  \over 2|z|  }\, |f(z)-f(0)|>1.$ 
D'apr\`es le th\'eor\`eme de Hahn-Banach, il existe une forme 
$ {\matC}$-lin\'eaire continue $ \lambda $ sur $ {\matC}^{n} $ de norme $ 1 $ telle
 que 
$$\lambda(f(0)+\zeta(f(z)-f(0))=\|f(0)+\zeta(f(z)-f(0))\|>1.$$
 L'application $ \lambda{\circ}f $ est bien holomorphe sur $ \Delta $ et
 $ (\lambda{\circ}f)(\Delta)\subset\overline{\Delta}$. On peut appliquer
 la premi\`ere partie de la d\'emonstration, ce qui donne la contradiction
 cherch\'ee. 
\medskip

On d\'efinit alors la notion de point complexe extr\'emal d'un
 ensemble $ A\subset{\matC}^{n}$. 
\medskip

{\bf D\'efinition 13.3. }- Soit $ A $ un partie de $ {\matC}^{n}$. On
 dit qu'un point $ x $ de $ A $ est un point complexe extremal de $ A $ si
 le seul vecteur $ y\in{\matC}^{n} $ tel que $ x+\zeta y $ appartienne
 \`a $ A$, pour tout $ \zeta\in{\matC},|\zeta|<1 $ est le vecteur nul.
 
\medskip

Sous les hypoth\`eses du lemme 13.2, si une application holomorphe
 $ f:\Delta\longrightarrow\overline{B} $ est telle que $ f(z)\neq
 f(0)$, alors l'image de $ \Delta $ par l'application lin\'eaire affine
 non constante 
$$\zeta\mapsto f(0)+\zeta\,{1-|z|  \over 2|z|
  }\, \hbox{\rm\ (}f(z)-f(0)\hbox{\rm)}$$
 est contenue dans $ \overline{B}$. Ceci montre la proposition suivante.
 
\medskip

{\bf Proposition 13.4. }{\it- Soit \/}$ f\in H(\Delta,{\matC}^{n})
 $ {\it telle que \/}$ f(\Delta)\subset\overline{B}. $ {\it Soit \/}$ z_{0}\in\Delta${\it.
 Si \/}$ f(z_{0}) $ {\it est un point complexe extr\'emal de \/}$\overline{B}${\it,
 alors \/}$ f $ {\it est constante. R\'eciproquement, si \/}$ z_{0}\in\overline{B}
 $ {\it n'est pas un point complexe extr\'emal de \/}$\overline{B}${\it,
 il existe une application holomorphe non constante \/}$ f\in H(\Delta,\overline{B})
 $ {\it telle que \/}$ f(0)=z_{0}.$ 
\medskip

Si on consid\`ere maintenant des applications d'un domaine $ D $ de
 $ {\matC}^{p} $ \`a valeurs dans $ {\matC}^{n}$, en coupant par des
 sous-espaces affines de dimension $ 1 $ et en utilisant le principe
 du prolongement analytique, on en d\'eduit le principe du maximum
 fort suivant (E. Thorp et R. Whitley [21]). 
\medskip

{\bf Th\'eor\`eme 13.5. }{\it- Soit \/}$ D $ {\it un domaine 
de \/}${\matC}^{p} $ {\it et soit \/}$ f\in H(D,{\matC}^{n}) $ {\it tel que \/}$ f(D)
 $ {\it soit contenu dans la boule-unit\'e ferm\'ee \/}$\overline{B}
 $ {\it de \/}${\matC}^{n} $ {\it(pour une norme \/}$\|.\|${\it). Supposons
 que tout vecteur de norme \/}$ 1 $ {\it est un point complexe extr\'emal
 de \/}$\overline{B}${\it. Alors \/}$ f(D)\subset B $ {\it ou \/}$ f $ {\it
 est constante.\/} 
\medskip

{\bf Exemple 13.6. }- Consid\'erons $ {\matC}^{n} $ muni de la norme
 hermitienne. Soit $ B $ sa boule-unit\'e ouverte. Alors les points complexes
 extr\'emaux de $ \overline{B} $ sont exactement les points de $ \partial
 B $ (c'est-\`a-dire les points de norme $ 1$). 
\medskip

{\it D\'emonstration. - \/}Ce r\'esultat se montre facilement par
 un calcul \'el\'ementaire. Si on n'aime pas les calculs \'el\'ementaires,
 on peut remarquer que, d'apr\`es le th\'eor\`eme de Krein-Milman,
 $ \overline{B}$, qui est compact est l'enveloppe convexe ferm\'ee de
 ses points extr\'emaux et a fortiori de ses points complexes extremaux.
 L'ensemble $ S $ de ses points complexes extr\'emaux est donc non vide
 et est contenu dans $ \partial B$. D'autre part, il est clair que $ S
 $ est stable sous l'action du groupe unitaire $ U(n) $ qui agit transitivement
 sur $ \partial B$. Ainsi $ \partial B=S.$ 
\medskip

{\bf14. G\'eod\'esiques complexes} 
\medskip

En suivant les id\'ees  de E. Vesentini [22], nous allons d\'efinir
 la notion de g\'eod\'esique complexe d'un domaine born\'e $ D $ de 
$ {\matC}^{n}$. 
\medskip

{\bf Th\'eor\`eme et d\'efinition 14.1. }{\it- Soit \/}
$\varphi:\Delta\longrightarrow D $ {\it une application holomorphe. Les conditions
 suivantes sont \'equivalentes :\/} 

(i) {\it il existe deux points distincts \/}$\zeta $ {\it et \/}$\eta
 $ {\it de \/}$\Delta $ {\it tels que \/} 

$c_{D}(\varphi(\zeta),\varphi(\eta))=c_{\Delta}(\zeta,\eta)=\omega(\zeta,\eta)
 ;$ 

(ii) {\it pour tous les points \/}$\zeta $ {\it et \/}$\eta $ {\it
 de \/}$\Delta${\it, on a :\/} 

$c_{D}(\varphi(\zeta),\varphi(\eta))=c_{\Delta}(\zeta,\eta)=\omega(\zeta,\eta)
 ;$ 

(iii) {\it il existe \/}$\zeta\in\Delta $ {\it et \/}$ v\in{\matC}
 (v\neq0) $ {\it tels que \/} 

$E_{D}(\varphi(\zeta),\varphi'(\zeta).v)=E_{\Delta}(\zeta,v)=|v|/(1-|\zeta|^{2})
 ;$ 

(iv){\it Pour tous \/}$\zeta\in\Delta $ {\it et \/}$ v\in{\matC} (v\neq0)${\it,
 on a : \/}$ E_{D}(\varphi(\zeta),\varphi'(\zeta).v)=E_{\Delta}(\zeta,v)=|v|/(1-|\zeta|^{2})${\it.\/}

{\it On dit alors que \/}$\varphi $ {\it est une g\'eod\'esique complexe
 de \/}$ D${\it.\/} 
\medskip

{\it Id\'ee de la d\'emonstration. - \/}Etant donn\'es deux points
 $ x $ et $ y $ de $ D $ (resp. $ x\in D $ et $ v\in{\matC})$, il existe une
 fonction holomorphe $ f:D\longrightarrow\Delta $ qui donne
 la distance de Carath\'eodory $ c_{D}(x,y) $ (resp. $ E_{D}(x,v)$). C'est
 une application simple du th\'eor\`eme de Montel. Supposons que (i)
 est v\'erifi\'e. On a donc : 
$$\omega(f(\varphi(\zeta),f(\varphi(\eta))=c_{D}(\varphi(\zeta),\varphi(\eta))=\omega(\zeta,\eta).$$
 D'apr\`es le lemme de Schwarz-Pick, $ f{\circ}\varphi $ est un automophisme
 analytique de $ \Delta$, et, quitte \`a composer $ f $ avec un automorphisme
 de $ \Delta $ bien choisi, on peut supposer que $ f{\circ}\varphi={\rm id}$.
 Ceci permet de montrer (ii), (iii) et (iv). Le reste de la d\'emonstration
 est laiss\'e en exercice. 
\medskip

La notion de g\'eod\'esique complexe est reli\'e \`a celle de r\'etraction
 holomorphe de la fa\c con suivante. 
\medskip

{\bf Th\'eor\`eme 14.2. }{\it- Soit \/}$ D $ {\it un domaine born\'e
 de \/}${\matC}^{n} $ {\it et soit \/}$ V $ {\it une sous-vari\'et\'e
 analytique ferm\'ee de \/}$ D $ {\it analytiquement isomorphe au disque-unit\'e
 \/}$\Delta${\it. Les propositions suivantes sont \'equivalentes :\/}

(i) $ V $ {\it est l'image d'une r\'etraction holomorphe \/}
$\rho:D\longrightarrow D $ {\it;\/} 

(ii) $ V $ {\it est l'image d'une g\'eod\'esique complexe \/}
$\varphi:\Delta\longrightarrow D $ {\it;\/} 

(iii) $ c_{D}|_{V}=c_{V} $ {\it;\/} 

(iv) {\it il existe deux points distincts \/}$ x $ {\it et \/}$ y $ {\it
 de \/}$ D $ {\it tels que \/}$ c_{D}(x,y)=c_{V}(x,y) $ {\it;\/} 

(v) $ E_{D}|_{T(V)}=E_{V} $ {\it;\/} 

(vi) {\it il existe \/}$ x\in D $ {\it et \/}$ v\in T_{x}(V), v\neq0${\it,
 tel que \/}$ E_{D}(x,v)=E_{V}(x,v)${\it.\/} 
\medskip

{\it D\'emonstration. - \/}Comme $ V $ est isomorphe au disque-unit\'e
 $ \Delta$, on sait d'apr\`es le th\'eor\`eme et d\'efinition 14.1 que
 (ii), (iii), (iv), (v) et (vi) sont \'equivalents. Montrons que (i)
 entra\^{\i}ne (iii). Soit $ i $ l'injection de $ V $ dans $ D $ et $ \rho{\circ}i
 $ est \'egal \`a ${\rm id}|_{V}$. On d\'eduit du fait que les distances invariantes
 sont contractantes les in\'egalit\'es suivantes : pour tous $ x $ et
 $ y $ appartenant \`a $ V$, 
$$c_{V}(x,y)\leq c_{D}(x,y)\leq c_{V}(x,y).$$
 On en d\'eduit l'\'egalit\'e annonc\'ee. 
\medskip

R\'eciproquement, montrons par exemple que (iv) entra\^{\i}ne (i).
 Soient $ x $ et $ y $ deux points distincts de $ V $ tels que $ c_{D}(x,y)=c_{V}(x,y)$.
 D'apr\`es le th\'eor\`eme de Montel, il existe une application holomorphe
 $ f:D\longrightarrow\Delta $ qui donne $ c_{D}(x,y)$, c'est-\`a-dire
 telle que 
$$c_{D}(x,y)=c_{\Delta}(f(x),f(y))=\omega(f(x),f(y)).$$
 Soit $ g $ un isomorphisme de $ \Delta $ sur $ V$. Alors 
$ g{\circ}f:D\longrightarrow V $ est telle que 
$$c_{V}(g(f(x)),g(f(y)))=c_{D}(x,y).$$
 Mais comme $ c_{V}(x,y)=c_{D}(x,y)$, on trouve que 
$$c_{V}(g(f(x)),g(f(y)))=c_{V}(x,y)\hbox{\rm.}$$
 D'apr\`es le lemme de Schwarz-Pick, $ g{\circ}f|_{V} $ est un automorphisme
 analytique de $ V$. Quitte \`a le composer avec un automorphisme analytique
 $ h $ de $ V $ bien choisi, on trouve 
$$h{\circ}g{\circ}f|_{V}=\hbox{\rm id}|_{V},$$
 $ h{\circ}g{\circ}f:D\longrightarrow D $ est bien une
 r\'etraction holomorphe de $ D $ sur $ V$. 
\medskip

Nous avons maintenant un th\'eor\`eme d'existence de g\'eod\'esiques
 complexes. 
\medskip

{\bf Th\'eor\`eme 14.3. }{\it- Soit \/}$ D $ {\it un domaine born\'e
 de \/}${\matC}^{n} $ {\it taut au sens de H. Wu \/}[24].{\it Pour qu'il existe
 une g\'eod\'esique complexe \/}$\varphi:\Delta\longrightarrow
 D $ {\it telle que deux points distincts \/}$ a $ {\it et \/}$ b $ {\it
 appartiennent \`a l'image \/}$\varphi(\Delta)${\it, il faut et il suffit
 que \/} 
$$c_{D}(a,b)=\delta_{D}(a,b).$$
 {\it D\'emonstration. - \/}Comme nous l'avons montr\'e au th\'eor\`eme
 14.2, si $ a $ et $ b $ appartiennent \`a l'image d'une g\'eod\'esique
 complexe $ \varphi:\Delta\longrightarrow D$, il existe
 une r\'etraction holomorphe $ \rho:D\longrightarrow\varphi(\Delta)$,
 et on a : 
$$\delta_{D}|_{\varphi(\Delta)}=\delta_{\varphi(\Delta)}, c_{D}|_{\varphi(\Delta)}=c_{\varphi(\Delta)}.$$
 Comme $ \varphi(\Delta) $ est isomorphe au disque-unit\'e $ \Delta $ et
 que $ c_{\Delta}=\delta_{\Delta}=\omega$, l'\'egalit\'e annonc\'ee
 s'en d\'eduit. 
\medskip

Montrons la r\'eciproque. Du fait que $ D $ est taut, on d\'eduit facilement
 que, si $ a $ et $ b $ sont deux points distincts 
de $ D$, il existe une fonction
 holomorphe $ \varphi:\Delta\longrightarrow D $ qui r\'ealise
 exactement $ \delta_{D}(a,b)$, c'est-\`a-dire qu'il existe $ \alpha
 $ et $ \beta $ appartenant \`a $ \Delta $ tels que $ \varphi(\alpha)=a,
 \varphi(\beta)=b $ et que 
$$\delta_{D}(a,b)=\delta_{D}(\varphi(\alpha),\varphi(\beta))=\omega(\alpha,\beta).$$
 Soit maintenant $ f:D\longrightarrow\Delta $ une application
 holomorphe telle que 
$$c_{D}(a,b)=\omega(f(a),f(b)).$$
 L'application $ f{\circ}\varphi $ de $ \Delta $ dans $ \Delta $ est telle
 que 
$$\omega(f(\varphi(\alpha)),f(\varphi(\beta)))=\omega(\alpha,\beta).$$
 C'est donc un automorphisme analytique du disque-unit\'e $ \Delta$,
 et $ \varphi:\Delta\longrightarrow D $ est une g\'eod\'esique
 complexe de $ D$. 
\medskip

En particulier, si $ D $ est un domaine born\'e convexe de $ {\matC}^{n}$,
 on sait d'ap\`es le th\'eor\`eme 5.9 que $ c_{D}=\delta_{D}$. On en
 d\'eduit qu'\'etant donn\'es deux points $ a $ et $ b $ de $ D$, il existe
 au moins une g\'eod\'esique complexe dont l'image contient les points
 $ a $ et $ b$. 
\medskip

{\bf15. Unicit\'e des g\'eod\'esiques complexes. Application aux points
 fixes d'applications holomorphes} 
\medskip

Comme nous l'avons d\'ej\`a dit, \'etant donn\'es deux points d'un
 domaine born\'e $ D$, il n'existe pas, en g\'en\'eral, de g\'eod\'esique
 complexe dont l'image contienne ces deux points. Un autre probl\`eme
 int\'eressant est celui de l'unicit\'e des g\'eod\'esiques complexes
 dont l'image contient deux points donn\'es. Pour cela, il nous faut
 d'abord remarquer que, si $ \varphi:\Delta\longrightarrow
 D $ est une g\'eod\'esique complexe de $ D $ et si $ f $ est un automorphisme
 analytque de $ \Delta$, $ \varphi{\circ}f $ est encore une g\'eod\'esique
 complexe de $ D$. Aussi, nous d\'efinirons l'unicit\'e des g\'eod\'esiques
 complexes de la fa\c con suivante : \'etant donn\'es deux points $ a
 $ et $ b $ appartenant \`a $ D $ (resp. $ a\in D $ et 
$ v\in T_{a}(D)={\matC}^{n}$), nous dirons qu'une g\'eod\'esique complexe 
$ \varphi:\Delta\longrightarrow D $ est l'unique g\'eod\'esique complexe passant
 par $ a $ et $ b $ (resp. passant par $ a $ et tangent \`a $ v$) si $ a $ et
 $ b $ appartiennent \`a $ \varphi(\Delta) $ (resp. $ a\in\varphi(\Delta)
 $ et $ v\in T_{a}(\varphi(\Delta))$) et si, pour tout autre g\'eod\'esique
 complexe $ \psi:\Delta\longrightarrow D $ v\'erifiant les
 m\^emes propri\'et\'es, il existe un automorphisme $ f $ de $ \Delta
 $ tel que $ \psi=\varphi{\circ}f$. 
\medskip

En particulier, L. Lempert [18] a montr\'e que, si $ D $ est un domaine
 born\'e strictement convexe \`a fronti\`ere suffisamment r\'eguli\`ere,
 alors les g\'eod\'esiques complexes existent et sont uniques. En fait,
 la r\'egularit\'e de la fronti\`ere n'est pas n\'ecessaire et S. Dineen
 ([7], p.93) montre le th\'eor\`eme suivant. On dit qu'un domaine $ D
 $ est strictement convexe s'il n'existe pas de segment $ [a,b] $ avec
 $ a\neq b $ contenu dans la fronti\`ere de $ D$. 
\medskip

{\bf Th\'eor\`eme 15.1. }{\it- Soit \/}$ D $ {\it un domaine born\'e
 strictement convexe de \/}${\matC}^{n}${\it. Alors les g\'eod\'esiques
 complexes de \/}$ D $ {\it sont uniques au sens pr\'ec\'edent.\/} 
\medskip

{\it D\'emonstration. - \/}Consid\'erons une g\'eod\'esique complexe
 $ \varphi:\Delta\longrightarrow D $ d'un domaine born\'e
 $ D $ de $ {\matC}^{n}$. Il est classique que $ \varphi $ admet des limites
 radiales 
$$\varphi^{*}(e^{i\theta})={\rm lim}_{r\rightarrow1_{-}}\varphi(re^{i\theta})$$
 presque partout. Supposons $ D $ strictement convexe. Faisons la d\'emonstration
 dans le cas de deux points $ a $ et $ b$. Soient $ \varphi $ et $ \psi $ 
deux g\'eod\'esiques complexes passant
 par $ a $ et $ b$. Quitte \`a les reparam\'etrer, on peut supposer qu'il
 existe $ \zeta $ et $ \eta $ tels que 
$$\varphi(\zeta)=\psi(\zeta)=a, \varphi(\eta)=\psi(\eta)=b.$$
 Soient $ \varphi^{*} $ et $ \psi^{*} $ les limites radiales d\'efinies
 presque partout. D'apr\`es le th\'eor\`eme 14.1, pour tout $ \lambda\in[0,1],
 \lambda\varphi+(1-\lambda)\psi $ est une g\'eos\'esique complexe de
 $ D $ passant par $ a $ et $ b$. Soit $ \theta\in{\matR} $ tel que $ \varphi^{*}(e^{i\theta})
 $ et $ \psi^{*}(e^{i\theta}) $ soient tous les deux d\'efinis. 
$ \lambda\mapsto\lambda\varphi^{*}(e^{i\theta})+(1-\lambda)\psi^{*}(e^{i\theta})
 $ est un segment contenu dans la fronti\`ere de $ D$. Comme $ D $ est
 strictement convexe, ceci entra\^{\i}ne que $ \varphi^{*}(e^{i\theta})=\psi^{*}(e^{i\theta})$,
 et ceci est vrai pour presque tout $ \theta$. Ceci entra\^{\i}ne que
 $ \varphi=\psi$. Le th\'eor\`eme est d\'emontr\'e. 
\medskip

Dans le cas de la boule-unit\'e de $ {\matC}^{n} $ pour une norme $ \|.\|$,
 on a le th\'eor\`eme d'unicit\'e suivant. 
\medskip

{\bf Th\'eor\`eme 15.2. }{\it- Soit \/}$ B $ {\it la boule-unit\'e ouverte
 de \/}${\matC}^{n} $ {\it pour une norme \/}$\|.\|${\it. Soit \/}$ v\in\partial
 B $ {\it un point complexe extr\'emal de \/}$\overline{B}${\it. Alors,\/}
 
$$\varphi:\zeta\mapsto\varphi(\zeta)=\zeta v$$
 {\it est l'unique g\'eod\'esique complexe de \/}$ B $ {\it telle que
 \/}$\varphi(0)=0, \varphi' (0)=v${\it.\/} 
\medskip

{\it D\'emonstration. - \/}Consid\'erons le d\'eveloppement en s\'erie
 de $ \varphi$ 
$$\varphi(\zeta)=\sum^{\infty}_{n=1}a_{n}\zeta^{n},$$
 avec $ a_{1}=v$. L'application $ f:\Delta\longrightarrow{\matC}^{n} $
 d\'efinie par $ f(\zeta)=\varphi(\zeta)/\zeta $ admet pour d\'eveloppement
 en s\'erie 
$$f(\zeta)=\sum^{\infty}_{n=1}a_{n}\zeta^{n-1}=v+\sum^{\infty}_{n=2}a_{n}\zeta^{n-1}.$$
 Pour toute forme lin\'eaire $ m $ sur $ {\matC}^{n} $ de norme $ 1$,
 $ m{\circ}\varphi $ s'annule en $ 0 $ et est \`a valeurs dans $ \overline{\Delta}$.
 Le lemme de Schwarz montre que $ m{\circ}\varphi $ est \`a valeurs dans
 $ \Delta $ et v\'erifie 
$$|m{\circ}\varphi(\zeta)|\leq|\zeta|.$$
 Par suite, $ m{\circ}f $ est \`a valeurs dans $ \overline{\Delta}$.
 Le th\'eor\`eme de Hahn-Banach montre que $ f $ envoie $ \Delta $ dans
 la boule-unit\'e ferm\'ee $ \overline{B} $ et $ f(0)=v$. D'apr\`es le
 th\'eor\`eme 13.5, $ f $ est constante \'egale \`a $ v$. Le th\'eor\`eme
 est d\'emontr\'e. 
\medskip

Nous allons maintenant appliquer ce r\'esultat aux points fixes d'applications
 holomorphes. 
\medskip

{\bf Th\'eor\`eme 15.3. }{\it- Soit \/}$ B $ {\it la boule-unit\'e ouverte
 de \/}${\matC}^{n} $ {\it pour une norme \/}$\|.\| $ {\it et supposons
 que tous les points de la fronti\`ere \/}$\partial B $ {\it de \/}$ B
 $ {\it dont des points complexe extr\'emaux de \/}$\overline{B}${\it.
 Soit \/}$ f:B\longrightarrow B $ {\it une application holomorphe
 telle que \/}$ f(0)=0${\it. Alors, l'ensemble des points fixes de \/}$ f
 $ {\it est l'intersection de \/}$ B $ {\it avec le sous-espace propre
 correspondant \`a la valeur propre \/}$ 1 $ {\it de \/}$ f' (0)${\it.\/}
 
\medskip

{\it D\'emonstration. - \/}On sait d\'ej\`a que  ${\rm Fix} f $ est une sous-vari\'et\'e
 analytique complexe de $ B $ et que son espace tangent en $ 0 $ est \'egal
 \`a $ E_{1}$, le sous-espace propre correspondant \`a la valeur propre
 $ 1 $ de $ f' (0)$. Soit $ v\in E_{1} $ de norme $ 1 $ et consid\'erons la
 g\'eod\'esique complexe $ \varphi:\Delta\longrightarrow
 B $ d\'efinie par $ \varphi(\zeta)=\zeta v$. Comme $ f' (0).v=v$, on a
 $ f{\circ}\varphi(0)=0, (f{\circ}\varphi)' (0)=v$. Par suite, $ f{\circ}\varphi
 $ est une g\'eod\'esique complexe de $ B$. D'apr\`es le th\'eor\`eme
 d'unicit\'e (th\'eor\`eme 15.2), ceci suffit \`a prouver que $ f{\circ}\varphi=\varphi$.
 On en \'eduit que  ${\rm Fix} f =B\cap E_{1}$. 
\medskip

Dans certains cas particuliers, on peut montrer un r\'esultat plus
 pr\'ecis (M. Herv\'e [11]). 
\medskip

{\bf Th\'eor\`eme 15.4. }{\it- Soit \/}$ B $ {\it la boule-unit\'e ouverte
 de \/}${\matC}^{n} $ {\it pour la norme hermitienne. Soit \/}
$ f:B\longrightarrow B $ {\it une application holomorphe. Si \/} ${\rm Fix} f $
 {\it est non vide, \/} ${\rm Fix} f ${\it est l'intersection de \/}$ B
 $ {\it avec un sous-espace affine de \/}${\matC}^{n}${\it.\/} 
\medskip

{\it Id\'ee de la d\'emonstration. - \/}On remarque d'abord que tous
 les points de la fronti\`ere $ \partial B $ de $ B $ sont des points complexes
 extremaux de $ \overline{B}$. Si on suppose de plus que $ f(0)=0$, alors
  ${\rm Fix} f $ est l'intersection de $ B $ avec un sous-espace vectoriel $ E
 $ de $ {\matC}^{n}$. Le cas g\'en\'eral se d\'emontre en utilisant
 le fait que $ B $ est homog\`ene. On en d\'eduit que  ${\rm Fix} f $ est l'image
 de $ B\cap E $ par un automorphisme analytique de $ B$. La forme explicite
 des automorphismes permet alors de montrer que l'ensemble des points
 fixes de $ f $ est l'intersection de $ B $ avec un sous-espace affine
 de $ {\matC}^{n}$. 
\medskip

{\bf Remarque 15.5. }- Dans le th\'eor\`eme 10.8, nous avons montr\'e
 que la boule-unit\'e $ B_{n} $ de $ {\matC}^{n} $ et le polydisque $ \Delta^{n}
 $ n'\'etaient pas analytiquement isomorphes. Pour cela, on \'etait
 ramen\'e \`a montrer qu'il n'existe pas d'isomorphisme lin\'eaire
 entre $ B_{n} $ et $ \Delta^{n}$. On peut aussi dire qu'un tel isomorphime
 n'existe pas parce que tous les points de la fronti\`ere de $ B_{n}
 $ sont des points complexes extr\'emaux alors que ce n'est pas le cas
 pour $ \Delta^{n}$. 
\medskip

{\bf16. Exemples et applications} 
\medskip

Soit $ D $ un domaine born\'e convexe de $ {\matC}^{n}$.
 Soit $ f:D\longrightarrow D $ une application holomorphe telle que
 ${\rm Fix} f $ soit non vide. Alors, si  ${\rm Fix} f $ est de dimension $ 0$,  ${\rm Fix} f $
 est r\'eduit \`a un point. Supposons maintenant que  ${\rm Fix} f $ est de
 dimension $ 1$. D'apr\`es les r\'esultats pr\'ec\'edents, il existe
 une r\'etraction holomorphe $ \rho:D\longrightarrow {{\rm Fix} f } $. En \'ecrivant les applications 
$$\hbox{\rm Fix }f\longrightarrow D\longrightarrow\hbox{\rm
 Fix }f\hbox{\rm,}$$
 et la suite d'homotopie associ\'ee, on montre que
le premier groupe d'homotopie $ \pi_{1}( {\rm Fix}  f)$ est nul.
 Ainsi,  ${\rm Fix} f $ est simplement connexe. Comme  ${\rm Fix} f $ admet des fonctions
 born\'ees non constantes,  ${\rm Fix} f $ est isomophe au disque-unit\'e $ \Delta$,
 c'est donc l'image d'une g\'eod\'esique complexe. 
\medskip

Consid\'erons maintenant le cas du bidisque $ \Delta^{2}$. On peut
 montrer facilement (E. Vesentini [22]) que, \`a un changement de param\`etre
 pr\`es, les g\'eod\'esiques complexes de $ \Delta^{2} $ sont d'une des
 deux formes suivantes 
$$\zeta\mapsto(\zeta,h(\zeta)) \hbox{\rm\ ou
 }\zeta\mapsto(h(\zeta),\zeta),$$
 o\`u $ h $ est une application holomorphe de de $ \Delta $ dans $ \Delta$.
        
\bigskip

\centerline{\bf Bibliographie} 
\medskip

1. M. Abate. Iteration theory, compactly divergent sequences and
 commuting holomorphic maps. Ann. Scuola Norm. Sup. Pisa Cl. Sci.
 (4), {\bf18 }(1991), 167--191. 
\medskip
2. T. Barth. The Kobayashi distance induces the standard topology. Proc. Amer. Math. Soc, {\bf35} (1972),
439--441.
\medskip

3. E. Bedford. On the automorphism group of a Stein manifold.
 Math. Ann., {\bf266 }(1983), 215--227. 
\medskip
4. N. Bourbaki, Int\'egration, chapitre 7, Hermann, Paris, 1963.
\medskip

5. H. Cartan, Les fonctions de deux variables complexes 
et le probl\`eme de la re\-pr\'esen\-tation analytique. 
J. Math. pures et appl., $9^e$ s\'erie, {\bf11} (1931), 1--114.
\medskip

6. H. Cartan. Sur les r\'etractions d'une vari\'et\'e. C.
 R. Acad. Sci. Paris Ser. I Math., {\bf303 }(1986), 715. 
\medskip

7. S. Dineen. The Schwarz lemma. Oxford Mathematical Monographs. Oxford Science Publications. The Clarendon Press,
Oxford University Press, New York, 1989.
\medskip
8. C. Earle and R. Hamilton. A fixed point theorem for holomorphic mappings. Proc. Symposia Pure Math., {\bf16} (1969),
 61--65.
\medskip

9. T. Franzoni and E. Vesentini. Holomorphic maps and invariant distances.
 Notas de Matematica [Mathematical Notes], {\bf69}. North-Holland Publishing
 Co, Amsterdam, 1980. 
\medskip
10. R.~Gunning and H.~Rossi, Analytic Functions of Several Complex Variables, Prentice-Hall, Inc., Englewood Cliffs, N.J. 1965
\medskip
11. M. Herv\'e. Quelques propri\'et\'es des applications analytiques des boules \`a $m$ dimensions dans elle-m\^eme.
J. Math. Pures et appl. (9), {\bf42} (1963), 117--147.
\medskip

12. L. Harris. Schwarz's lemma in normed linear spaces.
Proc. Nat. Acad. Sci. U.S.A. {\bf62} (1969), 1014--1017.

\medskip

13. L. Harris. Schwarz-Pick systems of pseudometrics for domains
 in normed linear spaces. Advances in holomorphy (Proc. Sem. Univ.
 Fed. Rio de Janeiro, Rio de Janeiro, 1977). North-Holland Math. Stud.,
 {\bf34}, 345-406. North-Holland, Amsterdam 1979. 
\medskip

14. M. Jarnicki and P. Pflug. Invariant distances and metrics in complex
 analysis. de Gruyter Expositions in Mathematics, {\bf9}, Walter de
 Gruyter Co, Berlin, 1993. 
\medskip

15. L. Kaup and B. Kaup. Holomorphic functions of several complex variables. An 
introduction to the fundamental theory. De Gruyter Studies in Math. {\bf3}, Berlin, 1983.
\medskip

16. S. Kobayashi. Hyperbolic complex spaces, Grundlehren 
der Mathematischen Wissenschaften [Fundamental Principles 
of Mathematical Sciences], 318, Springer-Verlag, Berlin, 1998.
\medskip

17. L. Lempert. La m\'etrique de Kobayashi et la 
repr\'esentation des domaines sur la boule.
 Bull. Soc. Math. France {\bf109} (1981), 427--474. 
\medskip

18. L. Lempert. Holomorphic retracts and intrinsic metrics in convex domains. 
Anal. Math. {\bf 8} (1982), 257--261. 
\medskip

19. P. Mazet and J.-P. Vigu\'e. Convexit\'e de la distance de Carath\'eodory et points fixes
 d'applications holomorphes. Bull. Sci. math. (2),
 {\bf116 }(1992), 285--305. 
\medskip
20. R. Narasimhan.  Several complex variables, The University of Chicago
 Press, Chicago, Ill., 1971. 
\medskip

21. E. Thorp and R. Whitley. The strong maximum modulus 
theorem for analytic functions into a Banach space.
 Proc. Amer. Math. Soc.  {\bf 18} (1967), 640--646.
\medskip
22. E. Vesentini. Complex geodesics and holomorphic maps.
 Symposia Mathematica, Vol. XXVI (Rome, 1980), 211--230. 
\medskip

23. J.-P. Vigu\'e. Sur les points fixes d'applications holomorphes. 
C. R. Acad. Sci. Paris S\'er. I Math. {\bf303} (1986),
 927--930.
\medskip

24. H.Wu, Normal families of holomorphic mappings, Acta.Math.{\bf 119}
1967, 193--233.
\medskip
\bigskip 

Jean-Pierre Vigu\'e 

UMR CNRS 6086 

Universit\'e de Poitiers 

Math\'ematiques 

SP2MI, BP 30179 

86962 FUTUROSCOPE  

e-mail : vigue@math.univ-poitiers.fr

ou jp.vigue@orange.fr

 \bye